\documentclass[a4paper,english,10pt,oneside]{article}\def\zibreport{1}\def\arxiv{1}

\usepackage[utf8]{inputenc}
\usepackage{hyperref}
\usepackage{amsmath, amsfonts}
\usepackage[pdftex,dvipsnames]{xcolor}
\usepackage{algorithm}
\usepackage{algpseudocodex}
\usepackage{xspace}
\usepackage{bm}
\usepackage{a4wide}
\usepackage{booktabs}
\usepackage{caption}
\usepackage{graphicx}
\usepackage{epstopdf}
\usepackage{tikz}
\usepackage{pgfplots}
\usepgfplotslibrary{statistics}
\pgfplotsset{compat=1.12}
\usepackage[export]{adjustbox}
\usepackage{cleveref}
\usepgfplotslibrary{external}
\usepackage{comment}
\usepackage{tablefootnote}
\usepackage{makeidx}
\usepackage{subcaption}

\usepackage[colorinlistoftodos,prependcaption,textsize=tiny]{todonotes}

\usepackage{orcidlink}

\newcommand{\myorcidlink}[1]{\,\href{https://orcid.org/#1}{\raisebox{-0.45ex}{\includegraphics[width=1.8ex]{orcid}}}}
\newcommand{\myurl}[1]{\textsf{\footnotesize \url{#1}}\xspace}

\newcommand{\allcaps}[1]{\protect\scalebox{0.93}{#1}}

\newcommand{\name}[1]{\mbox{#1}\xspace}
\newcommand{\nameCaps}[1]{\name{\allcaps{#1}}}

\newcommand{\remix}{\allcaps{REMix}\xspace}

\newcommand{\UNSEEN}{\nameCaps{UNSEEN}}
\newcommand{\SCIP}{\nameCaps{SCIP}}

\newcommand{\PySCIPOpt}{\nameCaps{PySCIPOpt}}

\newcommand{\CPLEX}{\nameCaps{CPLEX}}
\newcommand{\XPRESS}{\name{Xpress}}
\newcommand{\GUROBI}{\name{Gurobi}}

\newcommand{\COPT}{\nameCaps{COPT}}

\newcommand{\RAM}{\nameCaps{{RAM}}}
\newcommand{\GB}{\nameCaps{{GB}}}

\newcommand{\CPU}{\nameCaps{CPU}}

\newcommand{\LP}{\nameCaps{LP}}
\newcommand{\LPs}{\nameCaps{LPs}}
\newcommand{\MIP}{\nameCaps{MIP}}
\newcommand{\FP}{\nameCaps{FP}}

\newcommand{\MIPs}{\nameCaps{MIPs}}

\newcommand{\ESOM}{\nameCaps{ESOM}}
\newcommand{\ESOMs}{\nameCaps{ESOMs}}

\newcommand{\IPM}{\nameCaps{IPM}}
\newcommand{\IPMs}{\nameCaps{IPMs}}

\newcommand{\RENS}{\nameCaps{RENS}}

\newcommand{\GCNNs}{\nameCaps{GCNNs}}
\newcommand{\GCNN}{\nameCaps{GCNN}}

\newcommand{\NN}{\nameCaps{NN}}
\newcommand{\NNs}{\nameCaps{NNs}}
\newcommand{\NULL}{\nameCaps{NULL}}

\newcommand{\LPProb}{\ensuremath{\mathcal{P}_{\LP}(c, A, b, l, u)}}
\newcommand{\MIPProb}{\ensuremath{\mathcal{P}_{\MIP}(c, A, b, l, u, \mathcal{N}, \mathcal{I})}}

\newcommand{\T}{^{T}}

\newcommand{\subjectto}{\mbox{subject to}\xspace}

\newcommand{\floor}[1]{\left\lfloor{#1}\right\rfloor}
\newcommand{\ceil}[1]{\left\lceil{#1}\right\rceil}

\newcommand{\R}{\mathbb{R}}
\newcommand{\Rtimes}[2]{\R^{{#1}\times{#2}}}

\newcommand{\Rmn}{\Rtimes{m}{n}}

\newcommand{\Z}{\mathbb{Z}}

\definecolor{seagreen}{rgb}{0.18,0.74,0.56}
\definecolor{darkgreen}{rgb}{0.0,0.45,0.00}
\definecolor{navyblue}{rgb}{0.0,0.0,0.5}
\definecolor{steelblue}{rgb}{0.27,0.51,0.71}
\definecolor{siennabrown}{rgb}{0.63,0.32,0.18}
\definecolor{firebrickred}{rgb}{0.69,0.13,0.13}
\definecolor{gray75}{rgb}{0.75,0.75,0.75}
\definecolor{orange}{rgb}{.843,0.671,0.078}
\definecolor{gold}{rgb}{1.0,0.84,0.0}
\definecolor{scipyellow}{HTML}{FFFFD6}
\definecolor{soplexred}{HTML}{FFD8D8}
\definecolor{zimplgreen}{HTML}{D8FFD8}
\definecolor{ugblue}{HTML}{CFEFFF}
\definecolor{gcgorange}{HTML}{FFDAB9}

\definecolor{mDarkBrown}{HTML}{604c38}
\definecolor{mDarkTeal}{HTML}{23373b}
\definecolor{mLightBrown}{HTML}{EB811B}
\definecolor{mLightblue}{HTML}{14B03D}

\definecolor{c0}{HTML}{000060}
\definecolor{c1}{HTML}{0000FF}
\definecolor{c2}{HTML}{36648B}
\definecolor{c3}{HTML}{4682B4}
\definecolor{c4}{HTML}{5CACEE}
\definecolor{c5}{HTML}{FF0000}
\definecolor{c6}{HTML}{008888}
\definecolor{c7}{HTML}{00DD99}
\definecolor{c8}{HTML}{527B10}
\definecolor{c9}{HTML}{7BC618}
\definecolor{c10}{HTML}{8DD8F8}
\definecolor{background}{HTML}{FFFFFF}

\definecolor{cp1}{HTML}{F2AF29}
\definecolor{cp2}{HTML}{05668D}
\definecolor{cp3}{HTML}{02C39A}
\definecolor{cp4}{HTML}{4F345A}
\definecolor{cp5}{HTML}{F2B5D4}
\definecolor{cp6}{HTML}{DA2C38}

\renewcommand{\MIPProb}{\ensuremath{\mathcal{P}_{\MIP}(c, A, b, l, u, \mathcal{I})}}

\newcommand{\TheTitle}{Developing heuristic solution techniques for large-scale unit commitment models} 
\newcommand{\TheAuthors}{N.--C. Kempke, T. Kunt, B. Katamish, C. Vanaret, S. Sasanpour, J.--P. Clarner, T. Koch}

\hypersetup{%
    pdftitle={\TheTitle},
    pdfauthor={\TheAuthors}
}

\ifthenelse{\zibreport = 1}{
    \usepackage{./zibtitlepage}
}{}

\begin{document}

\ifthenelse{\zibreport = 1}{
    \ZTPTitle{\TheTitle}
    \title{\TheTitle}

    \ZTPAuthor{
        \ZTPHasOrcid{Nils--Christian Kempke}{0000-0003-4492-9818},
        \ZTPHasOrcid{Tim Kunt}{0009-0006-5732-3208},
        \ZTPHasOrcid{Bassel Katamish}{0009-0009-2345-9635}
        \ZTPHasOrcid{Charlie Vanaret}{0000-0002-1131-7631}
        \ZTPHasOrcid{Shima Sasanpour}{0000-0002-7502-6841}
        \ZTPHasOrcid{Jan--Patrick Clarner}{0000-0001-9412-5741}
        \ZTPHasOrcid{Thorsten Koch}{0000-0002-1967-0077}
    }
    \author{
        \ZTPHasOrcid{Nils--Christian Kempke}{0000-0003-4492-9818},\and\
        \ZTPHasOrcid{Tim Kunt}{0009-0006-5732-3208},\and\
        \ZTPHasOrcid{Bassel Katamish}{0009-0009-2345-9635},\and\
        \ZTPHasOrcid{Charlie Vanaret}{0000-0002-1131-7631},\and\
        \ZTPHasOrcid{Shima Sasanpour}{0000-0002-7502-6841},\and\
        \ZTPHasOrcid{Jan--Patrick Clarner}{0000-0001-9412-5741},\and\
        \ZTPHasOrcid{Thorsten Koch}{0000-0002-1967-0077}
    }

    \ZTPInfo{Preprint}
    \ZTPNumber{25-03}
    \ZTPMonth{February}
    \ZTPYear{2025}

    \date{\normalsize February 26, 2025}
    \ifthenelse{\arxiv = 0}{
        \zibtitlepage
    }{}
    \maketitle

}{}

\begin{abstract}
Shifting towards renewable energy sources and reducing carbon emissions necessitate sophisticated energy system planning, optimization, and extension.
Energy systems optimization models (\ESOMs) often form the basis for political and operational decision-making.
\ESOMs are frequently formulated as linear (\LPs) and mixed-integer linear (\MIP) problems.
\MIPs allow continuous and discrete decision variables.
Consequently, they are substantially more expressive than \LPs but also more challenging to solve.
The ever-growing size and complexity of \ESOMs take a toll on the computational time of state-of-the-art commercial solvers.
Indeed, for large-scale \ESOMs, solving the \LP relaxation -- the basis of modern \MIP solution algorithms -- can be very costly.
These time requirements can render \ESOM \MIPs impractical for real-world applications.
This article considers a set of large-scale decarbonization-focused unit commitment models with expansion decisions based on the \remix framework (up to 83 million variables and 900,000 discrete decision variables).
For these particular instances, the solution to the \LP relaxation and the \MIP optimum lie close.
Based on this observation, we investigate the application of relaxation-enforced neighborhood search (\RENS), machine learning guided rounding, and a fix-and-propagate (\FP) heuristic as a standalone solution method.
Our approach generated feasible solutions 20 to 100 times faster than \GUROBI, achieving comparable solution quality with primal-dual gaps as low as 1\% and up to 35\%.
This enabled us to solve numerous scenarios without lowering the quality of our models.
For some instances that \GUROBI could not solve within two days, our \FP method provided feasible solutions in under one hour.

\end{abstract}

{\centering\footnotesize In memory of Jan-Patrick Clarner.\par}

\section{Motivation}

The transition towards sustainable energy systems is one of the most pressing challenges of this century. As the global community strives to mitigate climate change and reduce greenhouse gas emissions, the mathematical development of energy system optimization models (\ESOMs) has become crucial.
Modeling various facets of energy systems in single analytical models should enable efficient planning and decision-making for adopting new policies.

Modern \ESOMs often capture as many aspects of an energy system as possible to increase the expressiveness of the model.
Power generation via various renewable sources such as solar, wind, and hydropower, as well as conventional sources such as fossil and nuclear power plants, are combined in an extensible power grid that spans vast regions \cite{WetzelGreenEnergyCarriers2023}.
Accurately modeling energy sources, network properties, and demand is fundamental for drawing reliable conclusions based on a given model \cite{DeCarolis_2017_BestPracticeModeling}.
The typical minimization of the total system costs leads to an affordable energy system. Additional constraints within the optimization problem, such as a carbon emission limit \cite{Victoria2020} and a minimum amount of backup capacities \cite{Sasanpour_2021}, can ensure that the energy system becomes sustainable and stays secure.
While the generation, transport and storage of power in \ESOM are usually modeled with linear approximations, \ESOMs with high spatial resolution should incorporate various discrete decisions \cite{Ommen_2014}.
These decisions can be split into several categories.
Operational decisions model the mode of the energy system, e.g., whether or not a power plant should be operated at a given point of discretized time (unit commitment) while considering the minimum load and minimum up- and down-times of the power plant \cite{HOFFMANN2024100190}.
Investment choices include extending the given energy systems, such as building new power plants or expanding the network with additional pipelines and grid lines.

Traditionally, \ESOMs are implemented as large-scale linear programming problems (\LPs) for which fast (in polynomial time) and reliable optimization algorithms exist.
However, they capture only continuous aspects of real-world scenarios and can only approximate discrete decisions.
On the other hand, mixed-integer linear programming problems (\MIPs) allow the modeling of discrete decisions, e.g., whether or not to expand or operate a specific energy source or pipeline. Their higher expressiveness enables modelers to describe real-world applications and connections more accurately.
This increase in eloquence when describing energy systems comes with the drawback of being substantially more difficult (NP-hard) to solve \cite{Wirtz_2021}.

Modern mixed-integer linear solvers such as \GUROBI \cite{Gurobi11}, \CPLEX \cite{CPLEX12}, \XPRESS \cite{FicoXpress9}, and \COPT \cite{COPT71} have been steadily improving at solving increasingly large \MIPs during the last decades \cite{kochProgressMathematicalProgramming2022}.
Still, in the face of the energy transition, \ESOMs have become increasingly complex.
Featuring decentralized girds with high spatial resolution, different energy sources, the varying availability of renewable energies., and complex decision-making, these models are often generated at the edge of what is still tractable for commercial software.
Thus, solution times play a vital role in the modeling process.

The \ESOMs we developed during the research project \UNSEEN\footnote{\UNSEEN project: \url{https://unseen-project.gitlab.io/home/}} follow a similar evolution: our large-scale \LP models evolved into \MIPs as our formulation developed over time.
However, at their full size, our latest models have become intractable for state-of-the-art \MIP solvers.
Fortunately, practitioners usually prefer meaningful, feasible solutions within an acceptable time since they can be produced in a fraction of the time required to find globally optimal solutions.
Moreover, models based on real-world predicted/historical data often contain uncertainty, which questions the relevance of globally optimal solutions.

In this article, we propose a problem-specific approach for solving large-scale \ESOMs as a competitive alternative to commercial optimization software: our heuristic techniques produce approximate, feasible solutions whose objective value is often close to the actual optimum of our models.
The (approximate) solution of multiple scenarios derived from the same model, obtained by sampling the input data over a given probability distribution, is used to mitigate uncertainty and improve robustness \cite{Frey2024_TacklingMultitudeOfUncertainties}.

The remainder of this article is structured as follows.
\Cref{sec:ESOM} briefly introduces the \MIP instances used in this article.
In \cref{sec:MIP}, we highlight the most important \MIP solving techniques and clarify concepts necessary for the following sections.
\Cref{sec:preliminary_analysis} discusses our preliminary analysis of the instances that motivated the development of tailored primal heuristics.
In \cref{sec:primal_heuristics}, we compare three primal heuristics designed to quickly find high-quality feasible solutions to the large-scale \MIPs.

\section{Modeling of energy systems}
\label{sec:ESOM}

The models considered in this article are based on the energy system optimization framework \remix \cite{REMIX_2012,REMIX_2017,Wetzel2024_REMIX}.
\remix is a feature-rich framework incorporating almost any temporal, spatial, and technological scale and detail.
Specifically, \remix enables the modeling of \ESOMs that can be used to analyze the decarbonization of our energy system.
The expansion and dispatch of different technologies, such as power plants, the energy grid, and energy storage, are typically optimized by minimizing the total system costs.
Usually, the energy system is modeled as an \LP.
However, if high spatial resolution is required, \MIP models that accurately describe strategies of individual power plants (such as their discrete expansion and unit commitment decisions for their operation) can be generated with \remix.

Our models' spatial basis is the German energy system. They describe the optimization of the power sector for 2030 with predefined capacities, e.g., for coal and lignite power plants, based on today's power plant park and the lifetime of individual power plants.
The energy system is optimized with hourly resolution, which results in 8,760 time steps.
A CO$_2$ price incentivizes the model to expand renewable energy technologies to reduce emissions.
Natural gas-fueled power plants are the only conventional technologies that can be further expanded apart from renewable energy sources.
Natural gas-fueled power plants and the grid can only be expanded discretely by adding individual power plants and transmission lines. The capacity of individual renewable energy power plants, such as solar and wind, is relatively small. Therefore, a continuous expansion is a sufficient approximation. To improve the representation of the operation of conventional power plants, a minimum partial load and minimum up- and downtime are considered.
A fixed time series from historical data represents the electricity exchange to Germany's neighboring countries.
The model consists of 488 nodes, representing the German power system at the transmission grid level (477 nodes) and its neighbors at the country level (11 nodes).
It can be aggregated on the spatial level, resulting in different sizes and difficulties. This is described in more detail in section \ref{sec:preliminary_analysis}.

Our \remix instances reflect the characteristics of \ESOMs: they incorporate multiple energy sources and address the expansion of and transition to renewable energies. Additionally, they have the advantage of offering different levels of complexity due to their scalability.

Since the input data, such as the weather data and the techno-economic parameters, are subject to uncertainties, a Monte Carlo approach was used to sample them.
Further information on the models, the input data, and the sampling can be found in \cite{dlr196232}.
\section{Mixed-integer linear programming background}
\label{sec:MIP}

An \LP with $n$ variables and $m$ linear constraints is written:
\begin{equation}
\label{eq:LinearProgram}
\tag{\LPProb}
\begin{alignedat}{2}
    \min \quad & c\T x \\
    \subjectto \quad & Ax \leq b \\
    & l \leq x \leq u,
\end{alignedat}
\end{equation}
where $x \in \R^n$ is a vector of real-valued variables, $(l, u) \in \R^n$ are the variable bounds, $A \in \Rmn$ the constraint matrix, $b\in\R^m$ the right-hand side and $c\in\R^n$ the objective gradient.
\LPs are well-studied optimization problems, and traditionally, the two most common solution techniques (leaving the rising family of first-order methods aside) are the simplex method \cite{DantzigWolfeSimplex_1955} and interior-point methods (\IPMs) \cite{terlaky2013_IPM, wright1997_primal_dual_ipm}.
While the worst-case complexity of the simplex method is exponential~\cite{KleeMintySimplexGood_1970}, on average, it performs in polynomial time~\cite{SpielmanTengSimplexPloyTime_2004}, which is the same complexity as \IPMs \cite{PotraWrightIPMs}.
Both methods can solve \LPs quickly and with high accuracy, and one might outperform the other on a given instance.

The extension from \LP to \MIP is seemingly simple. For a subset $\mathcal{I} \subset \{1,...,n\}$ of size $n_z \le n$, we additionally require that $x_\mathcal{I}$ be integral, where $x_\mathcal{I}$ is the sub-vector of $x$ consisting of all $x_i$ with $i\in \mathcal{I}$. A \MIP is thus written:
\begin{equation}
\label{eq:MixedIntegerLinearProgramm}
\tag{\MIPProb}
\begin{alignedat}{2}
    \min \quad & c\T x \\
    \subjectto \quad & Ax \leq b \\
    & l \leq x \leq u \\
    & x_\mathcal{I} \in \Z^{n_z}. 
\end{alignedat}
\end{equation}

The discrete (or integer) variables $x_\mathcal{I}$ are the combinatorial part of the optimization problem and allow the modeling of complex decisions (see \cref{sec:ESOM}).
We will often refer to $x_\mathcal{I}$ as $z$ to reduce the amount of indices (which also motivates the definition of $n_z$).

Classically, \MIPs are solved by a branch-and-bound approach \cite{LandDoigBnb_1960} that iteratively explores the combinatorial subspace of the optimization problem: an integer variable is picked at each level of the branch-and-bound tree and child nodes are generated, one for each possible integer assignment of the given variable.
This approach would yield an exponentially growing tree, effectively enumerating all possible integer solutions at its leaves.
In practice, lower and upper bounds on the globally optimal solution of the problem are maintained during the tree search, which allows certain portions of the tree to be pruned (e.g., if the lower bound of a sub-tree is larger than the best-known upper bound).
In particular, solving a relaxation of the \MIP problem at each node of the branch-and-bound tree yields a valid lower bound on the global optimal solution of the sub-tree originating at that node.
The standard \LP relaxation implemented in all commercial \MIP solvers is obtained by dropping the integrality constraints on the integer variables, which results in the LP \ref{eq:LinearProgram}.
On the other hand, upper bounds of the optimal solution can be obtained by primal heuristics \cite{MasterThesisTimo_2006} during the search: they are cheap strategies that determine satisfactory feasible points to the original problem, often by solving a subproblem of significantly reduced size.
A concise literature overview on the state of primal heuristics can be obtained via \cite{MasterThesisTimo_2006,PhDAchterberg,BertholdPhD2014,FischettiLodiHeurInMip2011}.

\begin{figure}[h!]
\centering
\includegraphics[width=10cm]{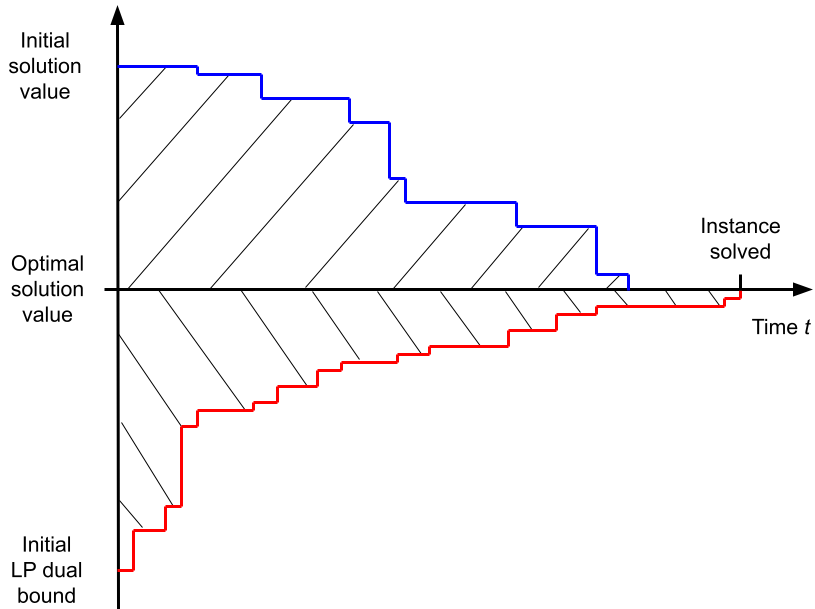}
\caption{Duality gap and primal-dual integral in mixed-integer programming.}
\label{fig:mip_gap}
\end{figure}

While solving \MIPs is NP-hard, the branch-and-bound approach is consistently being improved and often performs well in practice \cite{Achterberg2013, Lodi2010, bixby12brief}: the upper bounds (aka primal bounds) are enhanced by tailored primal heuristics.
Stronger problem formulations (e.g., with valid cuts) yield stronger lower bounds (also dual bounds).
The distance between upper and lower bounds, called \textit{duality gap}, monotonically decreases during the solution process (\cref{fig:mip_gap}).
In practice, users often set a positive \textit{target gap} (or gap tolerance) instead of solving the problem optimally with a zero gap but at a higher cost; the target gap measures the optimality of the best-known feasible point and limits the possible deviation to the optimal solution value.

The primal-dual integral \cite{BertholdPrimalDualIntegral_2013} describes the area between the curves of primal and dual bounds; it can be used to measure the impact of primal heuristics.
A smaller primal-dual integral generally means that good quality feasible points (given the available dual bound information) have been found more quickly.

The above discussion motivates the following nomenclature, standard in the \MIP community.
Solving a \MIP usually does not describe finding the (proven) global optimum but a feasible point within the target gap.
Feasible points found during the tree search are usually referred to as \textit{feasible solutions}, and any solutions satisfying the gap criterion are called (globally) \textit{optimal solutions}.
Usually, a \MIP solver will only produce a single optimal solution and multiple feasible solutions of varying quality during the solution process.
The objective value of the best optimal solution is called the \textit{optimum}.
The first feasible solution found during the optimization process is called \textit{initial solution}, and the currently best (with respect to its objective value) known feasible solution is called \textit{incumbent}.
We will use the same naming scheme in the remainder of this paper.
\section{Preliminary analysis of UNSEEN instances}
\label{sec:preliminary_analysis}

In this section, we derive some inherent characteristics of our models and their optimal solutions.
Commercial solvers struggled with or could not solve our most extensive energy system \MIPs, which made it hard to extract meaningful insight.
Instead, we generated smaller \ESOMs based on an aggregated underlying spatial resolution and varying input data for each scenario (see \cref{sec:ESOM}), which we could solve optimally with commercial solvers.
Our \MIP instances are listed in \cref{table:mipsizes}.
We refer to them as \textit{X-Small}, \textit{Small}, \textit{Medium}, \textit{Large}, and \textit{X-Large}, depending on the number of non-zeros in the constraint matrix $A$ in \ref{eq:MixedIntegerLinearProgramm}.
We generated 1,000 X-small and Small instances, 100 Medium and Large, and 20 X-Large instances.
X-Small, Small, and Medium instances are solved by commercial solvers within one hour.
Large instances are solved within hours; finding an initial solution for the X-Large instances may take multiple days.
Four Medium instances have proven infeasible; therefore, we are left with 96 for our analysis.

\begin{table}[htbp!]
\begin{center}
\begin{tabular}{c | c | c | c | c | c} 
& Variables & Integer variables & Constraints & Non-zeros & \#instances\\ 
\hline
X-Small & 67,386 & 5,445 & 73,518 & 244,293 & 1,000 \\ 
\hline
Small& 880,606 & 35,052 & 928,790 & 2,829,975 & 1,000 \\
\hline
Medium & 1,103,903 & 157,692 & 1,279,124 & 4,590,593 & 96 \\
\hline
Large & 13,597,290 & 972,477 & 16,313,022 & 47,310,667 & 100 \\
\hline
X-Large & 24,356,961 & 867,465 & 24,698,935 & 82,773,747 & 20 \\
\end{tabular}
\caption{Sizes of UNSEEN instances.}
\label{table:mipsizes}
\end{center}
\end{table}

We picked the 1,000 Small instances for our initial analysis, as these are the largest instances that can still be solved quickly.
For each Small instance, we used \GUROBI to solve the \MIP and the \LP relaxation (with the barrier method without crossover).
In the following, we will denote the optimal solutions to the \LP relaxation and to the \MIP by $x_{\text{\LP}}$ and $x_{\text{\MIP}}$, respectively.

\subsubsection*{The Initial Gap}

The optimum of the LP relaxation is given by $f^*_{\text{\LP}} := c\T x_{\text{\LP}}$.
Similarly, the optimum of the \MIP is given by $f^*_{\text{\MIP}} := c\T x_{\text{\MIP}}$.
We define the \textit{initial gap} $\Delta_{\text{init}}$ as the relative distance between the lower bound $f^*_{\text{\LP}}$ and the \MIP optimum $f^*_{\text{\MIP}}$, e.g.,
\begin{equation*}
    \Delta_{\text{init}} := \frac{\|f^*_{\text{\LP}} - f^*_{\text{\MIP}}\|}{\|f^*_{\text{\MIP}}\|}.
\end{equation*}
For $f^*_{\text{\LP}} =  f^*_{\text{\MIP}} = 0, \Delta_{\text{init}}$ is defined to be zero, for $f^*_{\text{\MIP}} = 0$ and $f^*_{\text{\LP}} \neq 0$,  $\Delta_{\text{init}}$ is defined to be $\infty$ (neither of these cases applied for any of our instances).
We use the initial gap as an intuition to guide our process: first, it gives an idea of how hard it might be to prove a feasible solution's optimality.
If the \MIP optimum and the initial dual bound are far apart, proving optimality might require extensive branching and cutting.
Second, if the \LP and \MIP optima are ``close'' (in the objective space), that is, when the initial gap is small, we hope this translates into the proximity of the \LP and \MIP optimal solutions.%
Note that computing the initial gap requires the knowledge of the \MIP optimum, which is the case only for instances of moderate size.

\begin{figure}[htbp!]
\centering
\begin{subfigure}{0.5\textwidth}
\centering
\tikzsetnextfilename{dualitygap_244k}
\begin{tikzpicture}
    \begin{axis}[
        ybar,
        width=0.95\textwidth,
        height=0.7\textwidth,
        ylabel = Frequency,
        xlabel = $\Delta_{\text{init}}$ in $\%$,
    ]
    \addplot+ [hist={bins=30}, fill=cp2!30, draw=cp2] table [y=duality_gap, col sep=comma] {data/duality_gap_in_percent_244k.txt};
    \end{axis}
\end{tikzpicture}
\end{subfigure}%
\begin{subfigure}[b]{0.5\textwidth}
\centering
\tikzsetnextfilename{dualitygap_2M}
\begin{tikzpicture}
    \begin{axis}[
        ybar,
        width=0.95\textwidth,
        height=0.7\textwidth,
        ylabel = Frequency,
        xlabel = $\Delta_{\text{init}}$ in $\%$,
    ]
        \addplot+ [hist={bins=30}, fill=cp2!30, draw=cp2] table [y=duality_gap, col sep=comma] {data/duality_gap_in_percent_2M.txt};
    \end{axis}
\end{tikzpicture}
\end{subfigure}
\caption{Distribution of \GUROBI's initial duality gap in percent for 1,000 X-Small (left) and Small (right) instances; average: $2.36 \times 10^{-6}$ (left) and $4.4 \times 10^{-2}$ (right).}
\label{fig:dualitygap_244_2M}
\end{figure}
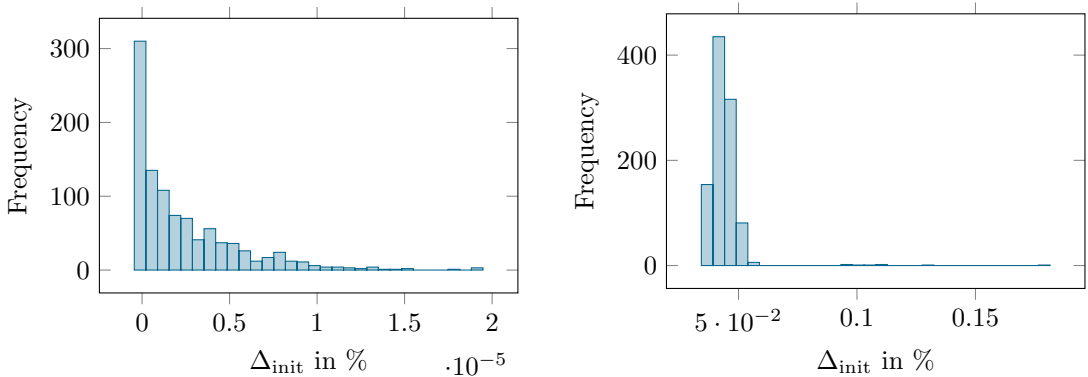

In \cref{fig:dualitygap_244_2M}, we plotted the frequency of initial gaps for the X-Small (left) and Small (right) instances solved with \GUROBI, with an average initial gap of $2.36 \times 10^{-6}$ and $4.4 \times 10^{-2}$, respectively.
For small instances, the relaxed \LP optimum and the \MIP optimum are likely to be close, although the initial gap seems to grow for larger instances.

\begin{figure}[h!]
\centering
\includegraphics[width=10cm]{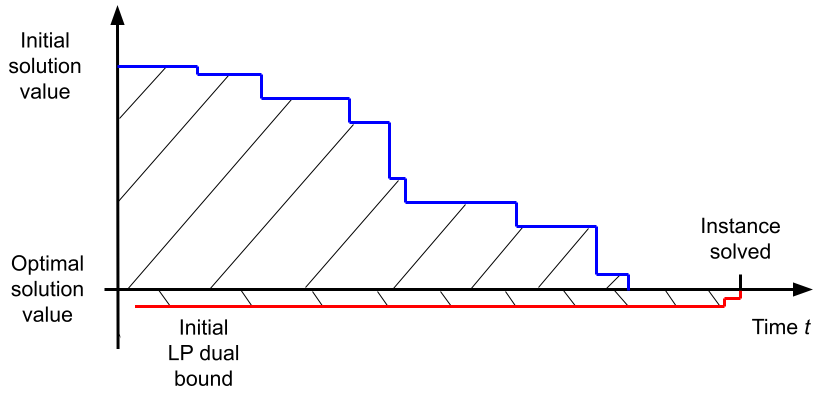}
\caption{Duality gap and primal-dual integral for \UNSEEN \MIP instances.}
\label{fig:unseen_mip_gap}
\end{figure}

The evolution of the duality gap for the UNSEEN instances is represented schematically in \cref{fig:unseen_mip_gap}.
In contrast to \cref{fig:mip_gap}, the primal-dual integral is mainly impacted by the primal bound.

These experiments highlight a particular characteristic of our models: the quality of the lower bounds computed by solving the \LP relaxation is such that branching and cutting should play little or no role in improving the dual bound.
Solving our \ESOMs mainly reduces to a purely primal problem: quickly finding good feasible points.

\subsubsection*{Solution distance}

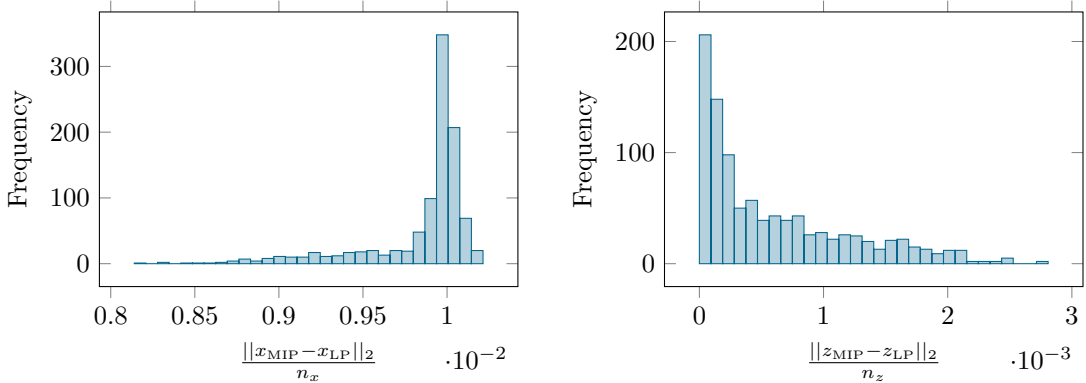
\begin{figure}[h!]
\centering
\begin{subfigure}{0.5\textwidth}
\centering
\tikzsetnextfilename{euclidean_distance_244k}
\begin{tikzpicture}
    \begin{axis}[
        ybar,
        width=0.95\textwidth,
        height=0.7\textwidth,
        ylabel = Frequency,
        xlabel = $\frac{||x_{\text{\MIP}}-x_{\text{\LP}}||_2}{n_x}$
    ]
        \addplot+ [hist={bins=30}, fill=cp2!30, draw=cp2] table [y=distance_all, col sep=comma] {data/euclidean_distance_244k.txt};
    \end{axis}
\end{tikzpicture}
\end{subfigure}%
\begin{subfigure}[b]{0.5\textwidth}
\centering
\tikzsetnextfilename{euclidean_integer_distance_244k}
\begin{tikzpicture}
    \begin{axis}[
        ybar,
        width=0.95\textwidth,
        height=0.7\textwidth,
        ylabel = Frequency,
        xlabel = $\frac{||z_{\text{\MIP}}-z_{\text{\LP}}||_2}{n_z}$
    ]
        \addplot+ [hist={bins=30}, fill=cp2!30, draw=cp2] table [y=distance_i, col sep=comma] {data/euclidean_distance_244k.txt};
    \end{axis}
\end{tikzpicture}
\end{subfigure}
\caption{Distribution of Euclidean distances over 1,000 X-Small instances; average: $9.85 \times 10^{-3}$ (left) and $5.98 \times 10^{-4}$ (right).}
\label{fig:euclidean_distance_244k}
\end{figure}

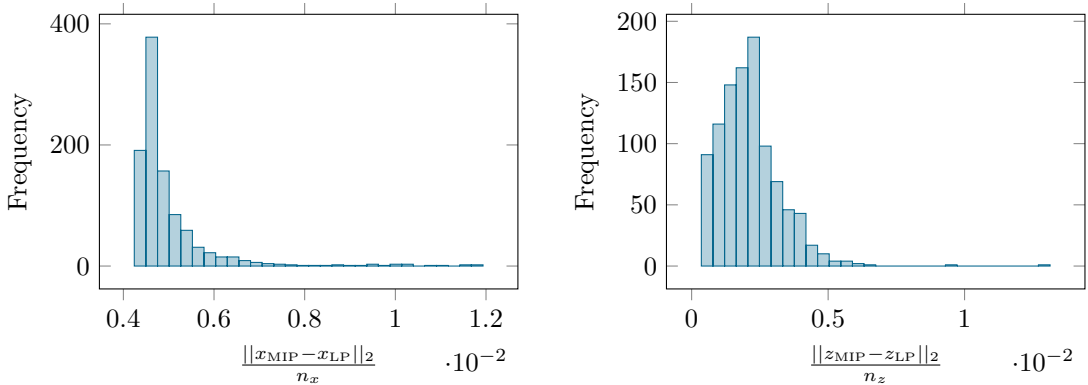
\begin{figure}[h!]
\centering
\begin{subfigure}{0.5\textwidth}
\centering
\tikzsetnextfilename{euclidean_distance_2M}
\begin{tikzpicture}
    \begin{axis}[
        ybar,
        width=0.95\textwidth,
        height=0.7\textwidth,
        ylabel = Frequency,
        xlabel = $\frac{||x_{\text{\MIP}}-x_{\text{\LP}}||_2}{n_x}$
    ]
        \addplot+ [hist={bins=30}, fill=cp2!30, draw=cp2] table [y=distance_all, col sep=comma] {data/euclidean_distance_2M.txt};
    \end{axis}
\end{tikzpicture}
\end{subfigure}%
\begin{subfigure}[b]{0.5\textwidth}
\centering
\tikzsetnextfilename{euclidean_integer_distance_2M}
\begin{tikzpicture}
    \begin{axis}[
        ybar,
        width=0.95\textwidth,
        height=0.7\textwidth,
        ylabel = Frequency,
        xlabel = $\frac{||z_{\text{\MIP}}-z_{\text{\LP}}||_2}{n_z}$
    ]
        \addplot+ [hist={bins=30}, fill=cp2!30, draw=cp2] table [y=distance_i, col sep=comma] {data/euclidean_distance_2M.txt};
    \end{axis}
\end{tikzpicture}
\end{subfigure}
\caption{Distribution of Euclidean distances over 1,000 Small instances; average: $4.99 \times 10^{-3}$ (left) and $2.12 \times 10^{-3}$ (right).}
\label{fig:euclidean_distance_2M}
\end{figure}

To assess the quality of the \LP relaxation solution of our \ESOM as a guide in finding near-optimal solutions, we investigate the distance between \LP relaxation solution and \MIP optimal solution.
\Cref{fig:euclidean_distance_244k} and \cref{fig:euclidean_distance_2M} show the average Euclidean distance between the \MIP optimal solution $x_{\text{\MIP}}$ and the \LP solution $x_{\text{\LP}}$ (left), and between the sub-vectors of integer variables $z_{\text{\MIP}}$ and $z_{\text{\LP}}$ (right) for X-Small and Small instances, respectively.
We observe a considerable proximity between \LP and \MIP solutions.
For the Small instances, the average distance between $x_{\text{\MIP}}$ and $x_{\text{\LP}}$ is $4.99 \times 10^{-3}$ and the average distance between $z_{\text{\MIP}}$ and $z_{\text{\LP}}$ is $2.12 \times 10^{-3}$.
This suggests that using the \LP relaxation solution as a starting point for a primal heuristic may be a sensible approach.

\begin{figure}[ht]
\centering
\begin{subfigure}{0.5\textwidth}
\centering
\tikzsetnextfilename{euclidean_distance_2M}
\begin{tikzpicture}
    \begin{axis}[
        ybar,
        width=0.95\textwidth,
        height=0.7\textwidth,
        ylabel = Frequency,
        xlabel = $\frac{||z_{\text{\MIP}}-z_{\text{\LP}}||_2}{n_z}$
    ]
        \addplot+ [hist={bins=30}, fill=cp2!30, draw=cp2] table [y=integral_fraction, col sep=comma] {data/integers_and_roundable_2M.txt};
    \end{axis}
\end{tikzpicture}
\end{subfigure}%
\begin{subfigure}[b]{0.5\textwidth}
\centering
\tikzsetnextfilename{euclidean_integer_distance_2M}
\begin{tikzpicture}
    \begin{axis}[
        ybar,
        width=0.95\textwidth,
        height=0.7\textwidth,
        ylabel = Frequency,
        xlabel = $\frac{||z_{\text{\MIP}}-\text{round}(z_{\text{\LP}})||_2}{n_z}$
    ]
        \addplot+ [hist={bins=30}, fill=cp2!30, draw=cp2] table [y=roundable_fraction, col sep=comma] {data/integers_and_roundable_2M.txt};
    \end{axis}
\end{tikzpicture}
\end{subfigure}
\caption{Distribution of percentages of coinciding integers before (left) and after rounding (right) over 1,000 Small instances; average: $82.39$ (left) and $94.21$ (right).}
\label{fig:coinciding_integers_2M}
\end{figure}
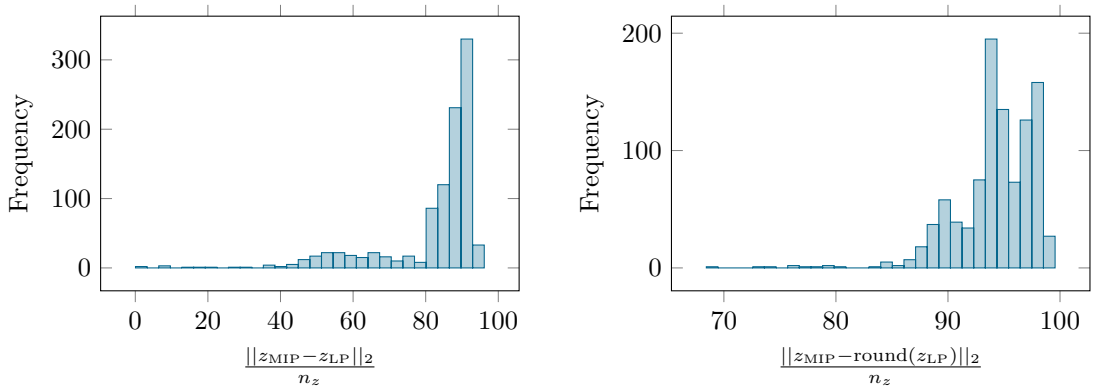

We analyze this relation even further and define as ``coinciding integers'' the integer variables that assume the same value in the \LP and \MIP optimal solutions.
Many coinciding integers would significantly increase the chance of finding high-quality \MIP feasible solutions from the starting point $x_{\text{\LP}}$.
As we generally cannot expect integer variables to lie strictly at integer values at a solution (this holds for both \LP and \MIP solvers), we consider that a variable in $x_{\text{\LP}}$ is equal to its \MIP counterpart if they differ by less than $10^{-6}$.
This is common practice in commercial/academic optimization software.
Usually, a \MIP solver considers a variable to be integer-feasible, e.g., to lie at an integer value, if its distance to the nearest integer is smaller than some threshold\footnote{For \GUROBI, this is $10^{-5}$, see the \GUROBI \href{https://docs.gurobi.com/projects/optimizer/en/current/reference/parameters.html\#parameterintfeastol}{manual}.}.
In \cref{fig:coinciding_integers_244k} and \cref{fig:coinciding_integers_2M}, we plotted the distribution of the percentage of coinciding integers before (left) and after (right) rounding each variable in $z_\text{\LP}$ to its nearest integer, for X-Small and Small instances, respectively.
For the X-Small instances, the average optimal solution consists of $93.39$ \% coinciding integers, while for the rounded solution, nearly all $99.76$ \% of integer variables coincide.
It is a notably high rate, but it is expected, given the earlier observations.
For the Small instances, the averages are still remarkably high, with $82.39$ \% before and $94.21$ \% after rounding.

\begin{figure}[ht]
\centering
\begin{subfigure}{0.5\textwidth}
\centering
\tikzsetnextfilename{integers_244k}
\begin{tikzpicture}
    \begin{axis}[
        ybar,
        width=0.95\textwidth,
        height=0.7\textwidth,
        ylabel = Frequency,
        xlabel = $\frac{||z_{\text{\MIP}}-z_{\text{\LP}}||_2}{n_z}$
    ]
        \addplot+ [hist={bins=30}, fill=cp2!30, draw=cp2] table [y=integral_fraction, col sep=comma] {data/integers_and_roundable_244k.txt};
    \end{axis}
\end{tikzpicture}
\end{subfigure}%
\begin{subfigure}[b]{0.5\textwidth}
\centering
\tikzsetnextfilename{roundable_244k}
\begin{tikzpicture}
    \begin{axis}[
        ybar,
        width=0.95\textwidth,
        height=0.7\textwidth,
        ylabel = Frequency,
        xlabel = $\frac{||z_{\text{\MIP}}-\text{round}(z_{\text{\LP}})||_2}{n_z}$
    ]
        \addplot+ [hist={bins=30}, fill=cp2!30, draw=cp2] table [y=roundable_fraction, col sep=comma] {data/integers_and_roundable_244k.txt};
    \end{axis}
\end{tikzpicture}
\end{subfigure}
\caption{Distribution of percentages of coinciding integers before (left) and after rounding (right) over 1,000 X-Small instances; average: $93.39$ (left) and $99.76$ (right).}
\label{fig:coinciding_integers_244k}
\end{figure}
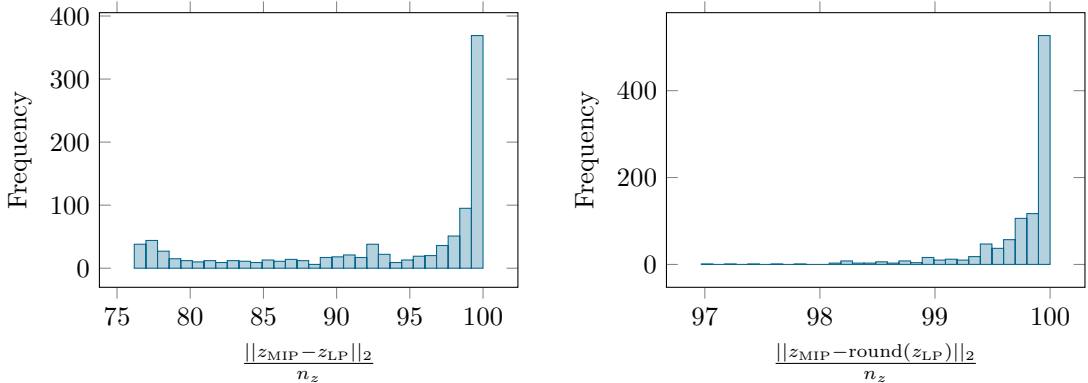

\newpage

\section{Primal heuristics for Energy System Optimization Models}
\label{sec:primal_heuristics}

In this section, we describe efficient primal heuristics for quickly finding good feasible points for our instances: i) a Relaxation Enforced Neighborhood Search heuristic, based on the proximity of the relaxed \LP and \MIP solutions (\cref{sec:RENS}); ii) a machine learning strategy for rounding the \LP solution (\cref{sec:learning-to-round}); and iii) an \LP-free fix-and-propagate strategy (\cref{sec:fix-and-propagate}).

\subsection{Relaxation enforced neighborhood search}
\label{sec:RENS}

Due to the proximity of the relaxed \LP solution and the \MIP optimal solution, high-quality feasible solutions (within the given tolerance) can be obtained by rounding fractional variables correctly.
Numerous \MIP heuristics exploited this (again refer to \cite{MasterThesisTimo_2006,PhDAchterberg,BertholdPhD2014,FischettiLodiHeurInMip2011}), including the Relaxation Enforced Neighborhood Search heuristic (\RENS) \cite{BertholdRENS2013}.
\RENS was a first stepping stone and benchmark toward the more elaborate primal heuristics discussed in the following subsections.

\subsubsection{Description of the heuristic}

Starting from the (fractional) \LP solution, \RENS \textit{enforces} a small neighborhood around each integer variable.
The resulting sub-\MIP has a drastically reduced search space.
Solving the \LP relaxation and the sub-\MIP may be significantly faster than solving the original \MIP while still producing near-optimal solutions.
A caveat is that incorrectly rounding or fixing a variable may lead to an infeasible sub-\MIP. 
Therefore, we restrict the integer variables to smaller neighborhoods around the fractional \LP solution, which keeps the sub-\MIP feasible.
We introduce additional constraints to the \MIP, a lower (resp. upper) bound by rounding down (resp. up) the \LP solution.
For non-fractional solution values, this directly leads to equality constraints.
The steps are summarized in \cref{alg:RENS}.
\begin{algorithm}[h!]
\caption{Relaxation enforced neighborhood search (\RENS)}\label{alg:RENS}
\begin{algorithmic}
\Require \MIPProb\ with integer variables $z := x_{\mathcal{I}}$
\Ensure \MIP solution or \NULL if sub-\MIP has no solution
\State Solve \LP relaxation and obtain optimal solution $x_{\text{\LP}}$
\State Impose constraints $\lfloor z_{\text{\LP}} \rfloor \leq z \leq \lceil z_{\text{\LP}} \rceil$\ on \MIP
\State Solve sub-\MIP
\If{sub-\MIP is feasible}
  \State \Return \MIP solution
\Else
  \State \Return \NULL
\EndIf
\end{algorithmic}
\end{algorithm}

Heuristics of this kind are part of every \MIP solver.
However, they tend to have much overhead and are usually employed with other heuristics.
Furthermore, successful branch-and-bound \MIP solvers usually require a crossover after the root node.

\subsubsection{Numerical results}

We now compare \RENS (using \GUROBI as \LP and \MIP solver) against standalone \GUROBI.
Computational time is improved if solving the \LP relaxation and the sub-\MIP is faster than solving the original \MIP.
All Small, Medium, Large, and X-Large instances were solved on a cluster of PowerEdge R650 machines, running each instance on a single Intel Xeon Gold 6342 \CPU with 2.8 GHz and 200 \GB \RAM.

Within the time limit, \RENS solved all Small and Medium instances and $88\%$ of the Large instances.
\GUROBI solved all Small and Medium instances, and $72\%$ of the Large instances (see \cref{table:RENS-results-multiple}).
On the X-Large set, \GUROBI could not terminate with an optimal \MIP solution for any instance.
After rounding, no \RENS sub-\MIP was detected as infeasible for any instance, a hint that the neighborhood could be tightened more aggressively.
The \RENS optima all lie within a margin of $0.0001\%$ to the \MIP optima (well within \GUROBI's default \MIP gap).
The quality of the \RENS solutions is identical to \GUROBI's on the original \MIP.

\begin{table}[h!]
\small
\begin{center}
\begin{tabular}{llccrr}
\toprule
Instance set & Solver & Solved instances & Success rate & Mean time [s] & Mean time [s] \\
             &        &                  & (\%) & (arithmetic) & (geometric) \\
\midrule
Small    & \RENS     & 1,000 & 100 &   38.57 &  37.59 \\
         & \GUROBI & 1,000 & 100 &   34.91 &  28.83 \\
\midrule
Medium   & \RENS     &   96 & 100 &   141.17 & 135.85 \\
         & \GUROBI &   96 & 100 &   867.38 &  764.55 \\
\midrule
Large    & \RENS     &  88 & 88 &  28,708.02 &  25,666.24 \\
         & \GUROBI &   72 &  72 & 23,924.49 &  20,315.41 \\
\bottomrule
\end{tabular}
\caption{RENS numerical results on multiple instance sets.}
\label{table:RENS-results-multiple}
\end{center}
\end{table}

\begin{figure}[ht]
    \centering
    \begin{subfigure}{0.32\textwidth}
        \centering
        \includegraphics[width=\linewidth]{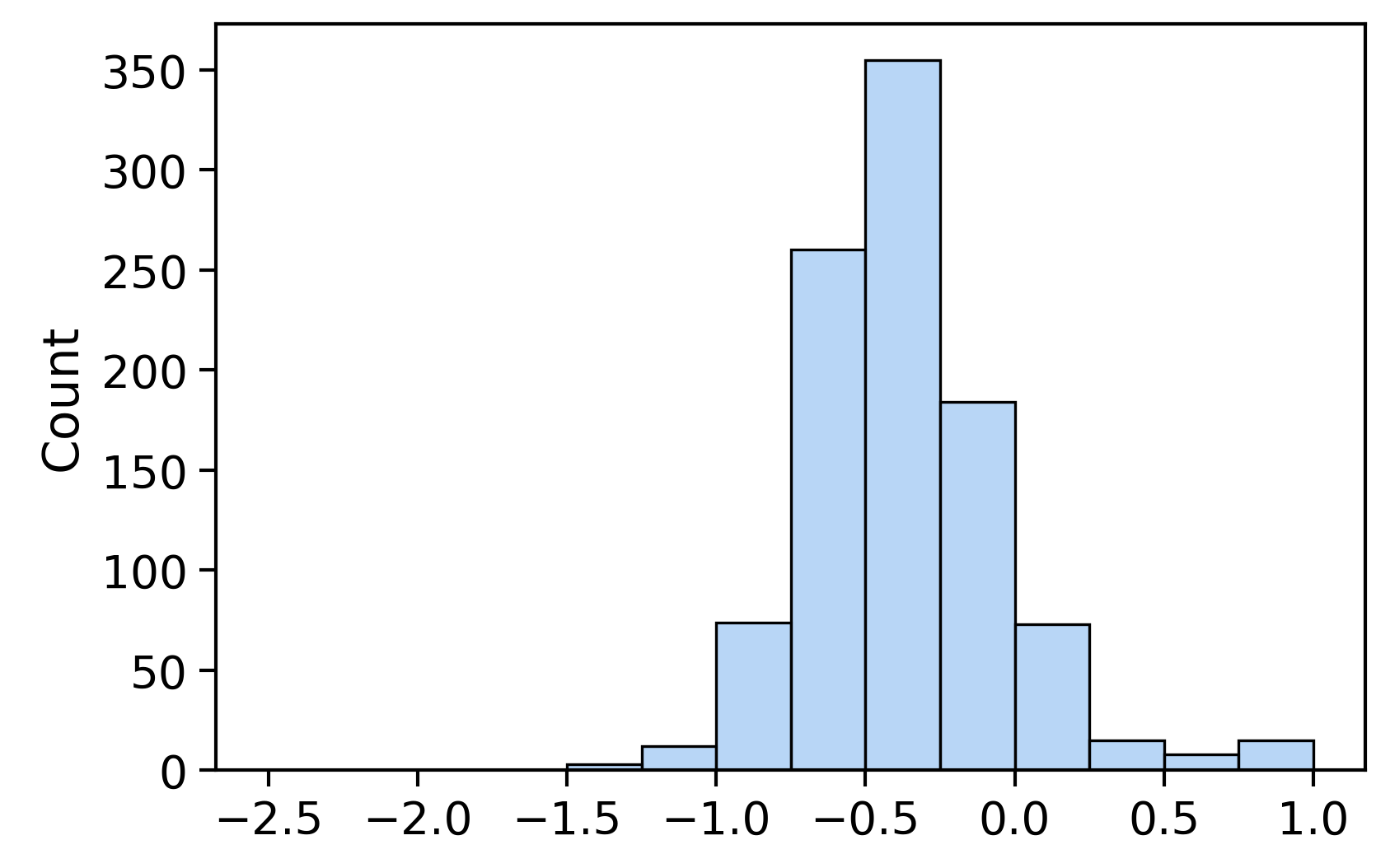}
        \caption{Small}
        \label{fig:rensimage3}
    \end{subfigure}
    \hfill
        \begin{subfigure}{0.32\textwidth}
        \centering
        \includegraphics[width=\linewidth]{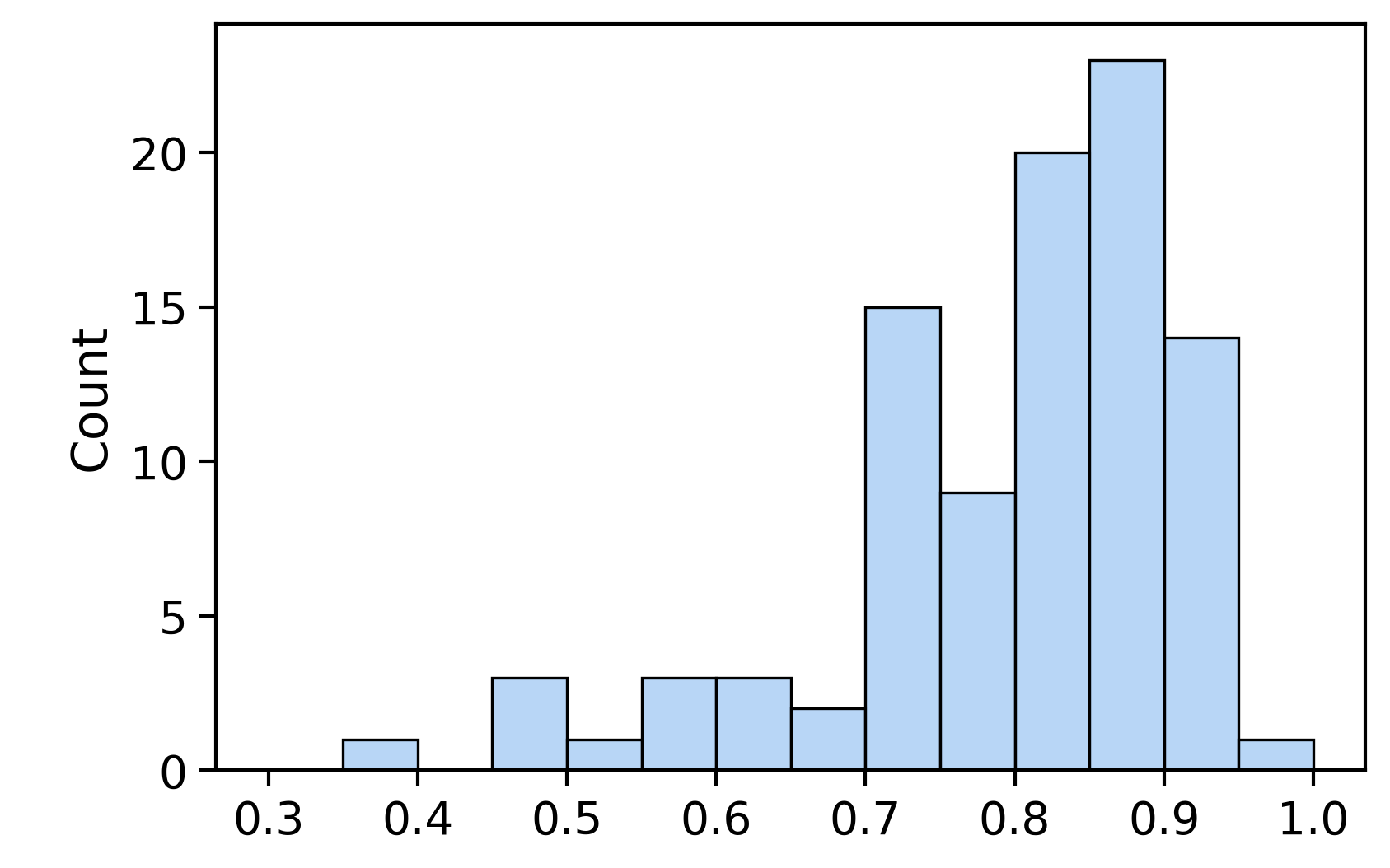}
        \caption{Medium}
        \label{fig:rensimage1}
    \end{subfigure}
    \hfill
    \begin{subfigure}{0.32\textwidth}
        \centering
        \includegraphics[width=\linewidth]{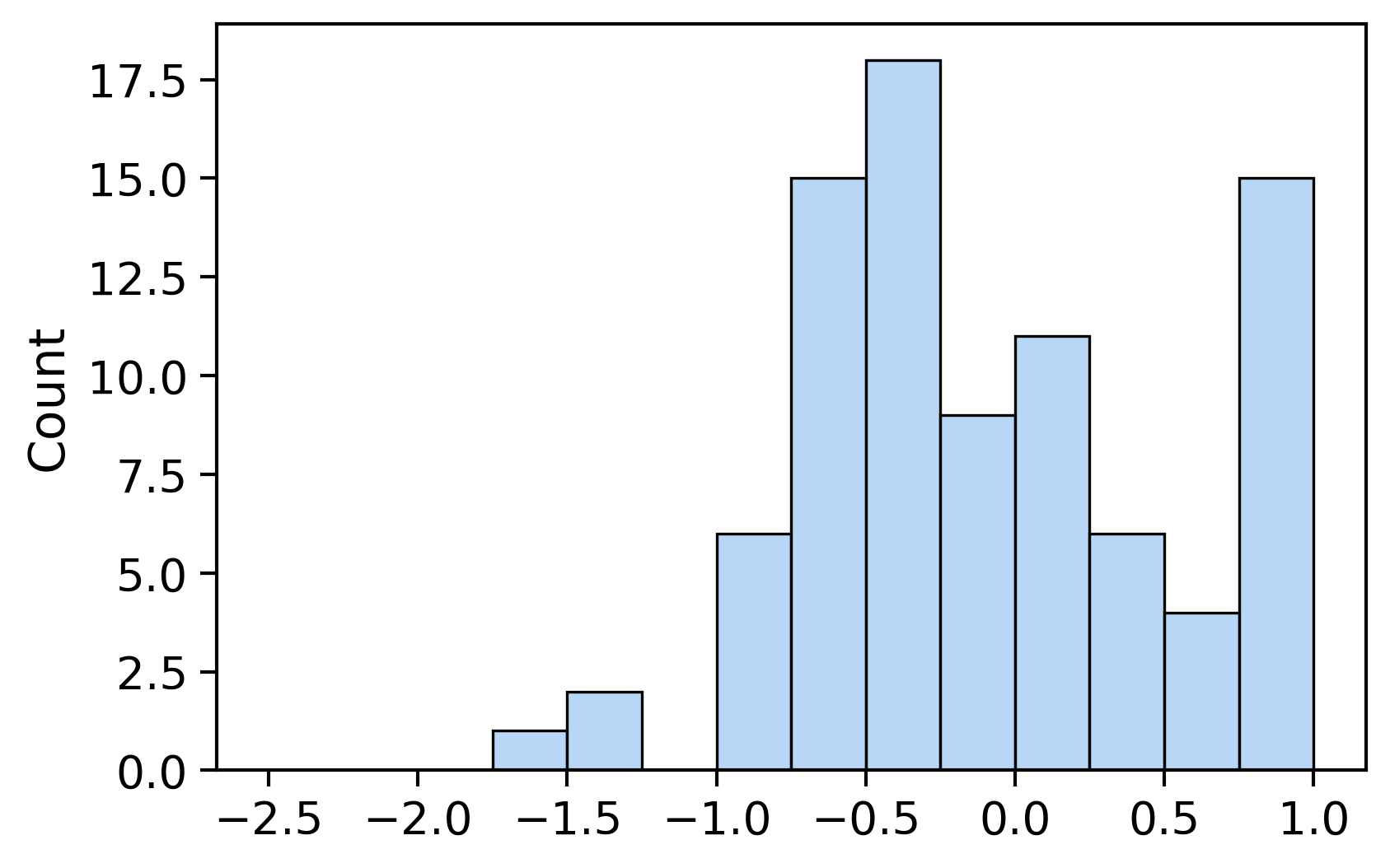}
        \caption{Large}
        \label{fig:rensimage4}
    \end{subfigure}
    \caption{relative time save of \RENS compared to \GUROBI on different instance sets.}
    \label{fig:RENS-time-save}
\end{figure}

The relative time-save of \RENS compared to \GUROBI for all three instance sets varies greatly, as shown in figure \ref{fig:RENS-time-save}
For the Small instances, the solution time for the \LP relaxation is often not significantly shorter than for the \MIP itself; consequently, \RENS runs longer than \GUROBI.
\RENS shows its potential on Medium instances with an overall average speedup of $79.20\%$ (the minimum being $29.55\%$ and the maximum $95.59\%$). On Large instances, \RENS gets mixed results; however, it solves more instances than \GUROBI, breaking additional instances that \GUROBI was not able to solve within the time limit.
For the unsolved instances, the \LP relaxation was solved within the time limit with a maximum of 88,739s and a mean of 66,565s, but the sub-\MIP could not be solved within the remaining time for all instances.

In this section, we demonstrated the promising potential of primal heuristics that exploit the rounding or fixing of integer variables. \RENS exhibited significant speedup on the Medium instances, higher robustness than \GUROBI on Large instances, and could generally produce high-quality solutions. However, the necessity to solve both the \LP relaxation and the sub-\MIP limits its competitiveness.

\subsection{Learning to round: a machine learning heuristic}
\label{sec:learning-to-round}

Since each instance set represents different scenarios for the same spatial aggregation, the structure of the instances -- such as the number of integers -- is identical within each instance set.
Therefore, it is reasonable to assume that machine learning-guided heuristics may perform well at identifying suitable rounding candidates and directions.

In this section, we describe a machine learning-guided heuristic that learns, for a given instance, whether or not a particular variable of the relaxed \LP solution should be rounded up or down to obtain an optimal \MIP solution.
We trained a multi-layer perceptron neural network (\NN) whose input is $z_\text{\LP}$, the key numerical values influencing the rounding decision. 
The output is a classification of variables, where +1 corresponds to ``round up'', -1 to ``round down'', and 0 to ``leave alone'' for variables that cannot be rounded.
The output layer uses Tanh units with $[-1, 1]$ range.
The \NN architecture, including the number of layers and hidden units per layer, was chosen after a validation phase: larger \NNs could not capture more information and did not generalize better.
To balance the architecture complexity and the generalization capabilities, we picked the smallest \NN that still achieved good generalization on the validation set while not overfitting on the training data.
We opted for three hidden layers with 300/400 ReLU units each; this enables the \NN to capture non-linear patterns efficiently and avoids issues like vanishing gradients that can arise with \NNs.

We restricted ourselves to the three instance sets X-Small, Small, and Large, and generated an additional 900 Large instances for training.
Each set of 1,000 instances was split into 700 instances for training, 200 for validation, and 100 for evaluation.

\subsubsection{Training of the neural network}

The training data was generated by solving all instances with \GUROBI twice: once as an \LP and once as a \MIP. 
Next, we extracted the integer parts of each solution, compared the two vectors, and generated target vectors that correspond to rounding instructions and serve as labels for the training phase.
The \NN was trained independently for each instance set.
We used an $\ell_1$ loss function and the Adam optimizer~\cite{Kingma2014}, implemented using the machine learning library PyTorch~\cite{PyTorch}.
The training results are summarized in \cref{table:mlresults}.
A classification accuracy of 96.1\% means that, for a given X-Small instance, the \NN should, on average, be able to correctly round 96.1\% of the variables.

\begin{table}[h!]
\footnotesize
\begin{center}
\begin{tabular}{c | c | c} 
& (Hidden Layers) $\times$ (Neurons) & Classification accuracy (\%) \\ 
\hline
X-Small & 3 $\times$ 300 & 96.1 \\ 
\hline
Small   & 3 $\times$ 300 & 81.0 \\
\hline
Large   & 3 $\times$ 400 & 82.7 \\
\end{tabular}
\caption{Machine learning heuristic: results of the training phase.}
\label{table:mlresults}
\end{center}
\end{table}

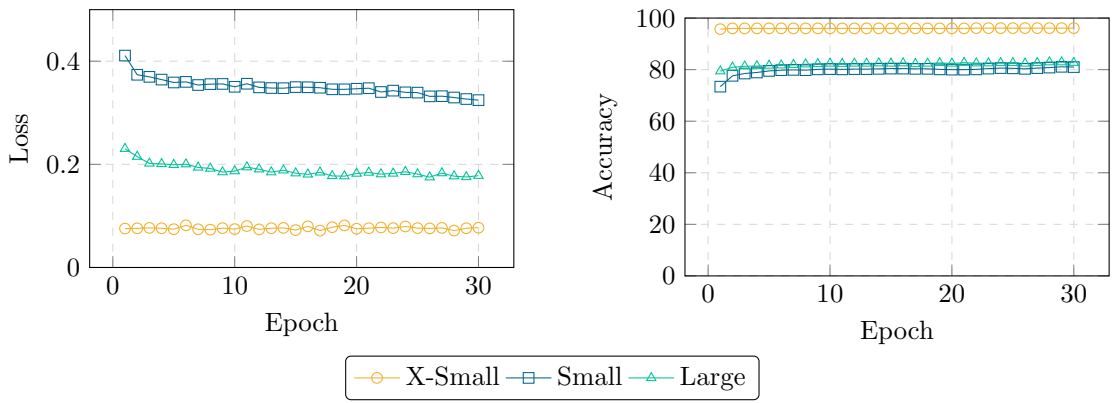
\begin{figure}[htbp!]
\centering
\begin{minipage}{0.48\textwidth}
    \centering
    \begin{tikzpicture}
        \begin{axis}[
            width=\textwidth, height=5cm,
            xlabel={Epoch}, ylabel={Loss},
            ymin=0,
            ymax=0.5,
            legend style={at={(0.5,-0.2)},anchor=north,legend columns=-1},
            grid=both,
            grid style={dashed, gray!30},
            legend style={draw=none}, 
        ]
            \addplot[color=cp1, mark=o] table[x=epoch, y=loss, col sep=comma] {data/ml_training_small.csv};
            \addplot[color=cp2, mark=square] table[x=epoch, y=loss, col sep=comma] {data/ml_training_mid.csv};
            \addplot[color=cp3, mark=triangle] table[x=epoch, y=loss, col sep=comma] {data/ml_training_large.csv};
        \end{axis}
    \end{tikzpicture}
\end{minipage}\hfill
\begin{minipage}{0.48\textwidth}
    \centering
    \begin{tikzpicture}
        \begin{axis}[
            width=\textwidth, height=5cm,
            xlabel={Epoch}, ylabel={Accuracy},
            ymin=0,
            ymax=100,
            legend style={at={(0.5,-0.3)},anchor=north,legend columns=-1},
            grid=both,
            grid style={dashed, gray!30},
            legend style={draw=none}, 
        ]
            \addplot[color=cp1, mark=o] table[x=epoch, y=accuracy, col sep=comma] {data/ml_training_small.csv};
            \addplot[color=cp2, mark=square] table[x=epoch, y=accuracy, col sep=comma] {data/ml_training_mid.csv};
            \addplot[color=cp3, mark=triangle] table[x=epoch, y=accuracy, col sep=comma] {data/ml_training_large.csv};
        \end{axis}
    \end{tikzpicture}
\end{minipage}
\begin{tikzpicture}
    \begin{axis}[
        hide axis,
        xmin=0, xmax=1, ymin=0, ymax=1,
        legend columns=3,
        legend style={at={(0.5,-0.2)}, anchor=north, draw=black, 
                      rounded corners=2pt}, 
    ]
       \addlegendimage{cp1, mark=o}
        \addlegendentry{X-Small}

        \addlegendimage{cp2, mark=square}
        \addlegendentry{Small}

        \addlegendimage{cp3, mark=triangle}
        \addlegendentry{Large}
    \end{axis}
\end{tikzpicture}

\caption{Loss and accuracy vs. epochs during training for X-Small, Small, and Large sets.}
\label{fig:loss_accuracy_combined}
\end{figure}

\Cref{fig:loss_accuracy_combined} depicts the training loss and the accuracy over the epochs.
Our experiments indicate that the classification error increases as the instance size increases but does not seem to be proportional to it.
This raises hope for the method's applicability.

\subsubsection{Description of the heuristic}

Although highly accurate, the rounding instructions produced by the \NN may fix integer variables to incorrect values, which causes the sub-\MIP to be infeasible.
This is a common issue for rounding heuristics, which we addressed by embedding our \NN into an algorithm using a fix-and-propagate scheme (\cref{alg:machine-learning-fix-and-propagate}).

\begin{algorithm}[h!]
\caption{Integer fixing using machine learning.}
\label{alg:machine-learning-fix-and-propagate}
\begin{algorithmic}
\Require \MIP $\MIPProb$, pre-trained neural network $\Theta$
\Ensure \MIP feasible solution or NULL if no solution was found
\State Presolve \MIP using \SCIP
\State Solve \LP relaxation: $(x_{\text{\LP}}, z_{\text{\LP}}) \gets \text{solve } \LPProb$
\State Compute rounding instructions: $\phi \gets \Theta(z_{\text{\LP}})$
\State $(\overline{l}, \overline{u}) \gets (l, u)$
\For{$i \in \mathcal{I}$}
  \If{$l_i = u_i$}
    \textbf{continue}
  \EndIf
  \State Round the bounds: 
  \begin{equation*}
  (\overline{l}_i, \overline{u}_i) \gets
  \begin{cases}
  (\ceil{{x_{\text{\LP}}}}, \ceil{{x_{\text{\LP}}}}) & \text{if } \phi_i = 1 \\
  (\floor{{x_{\text{\LP}}}}, \floor{{x_{\text{\LP}}}}) & \text{if } \phi_i = -1 \\
  (l_i, u_i) & \text{if } \phi_i = 0
  \end{cases}
  \end{equation*}
  \State $(\text{status}, \overline{l}, \overline{u}) \gets \text{propagate}(A, b, l, u,\mathcal{I}, \overline{l}, \overline{u})$ \Comment{\Cref{alg:propagation}.}
  \If{status is INFEASIBLE}
    \State \textbf{continue}
  \Else
    \State $(l, u) \gets (\overline{l}, \overline{u})$
  \EndIf
\EndFor
\If{$l_i = u_i, \forall i \in \mathcal{I}$}
  \State Solve resulting \LP to obtain a \MIP solution: $x_{\text{\MIP}} \gets \text{solve } \LPProb$
  \State Postsolve $x_{\text{\MIP}}$ with SCIP
  \State \Return  $x_{\text{\MIP}}$ 
\Else
  \State \Return \NULL
\EndIf
\end{algorithmic}
\end{algorithm}

After the \LP relaxation is solved, the \NN is fed the integer part of the \LP solution and produces the rounding instructions.
Then, we fix each unfixed integer variable by rounding it according to the instructions.
After each fixing, we propagate the bound change (\cref{alg:propagation}) to the rest of the bounds.
If the propagation phase detects an infeasible variable fixing, the fixing is dropped, and we move on to the next variable.
When all integer variables have been fixed (actively through rounding or passively during propagation), we solve the resulting \LP to obtain the optimal values for the continuous variables.

\begin{algorithm}[h!]
\caption{Propagation.}
\label{alg:propagation}
\begin{algorithmic}
\Require $\MIPProb$ and tightened bounds $l \leq \overline{l}$, $\overline{u} \leq u$.
\Ensure Whether $\ensuremath{\mathcal{P}_{\MIP}(c, A, b, \overline{l}, \overline{u},\mathcal{I})}$ is infeasible and updated bounds $\overline{l}$, $\overline{u}$.
\State Determined updated variables $Q = \{j : l_j \neq{} \overline{l}_j \lor u_j \neq \overline{u}_j\}$
\While {$Q \neq \emptyset$}
  \State Pop $j$ from $Q$
  \State Collect affected rows $R = \{i : A_{ij} \neq 0\}$
  \For{$i \in R$} \Comment{Propagate each affected row.}
    \State $\text{act}_- = \min\limits_{\overline{l}\leq x\leq \overline{u}} \{\sum\limits_{A_{ij} > 0} A_{ij} x_j + \sum\limits_{A_{ij} < 0} A_{ij} x_j\} = \sum\limits_{A_{ij} > 0} A_{ij} \overline{l}_j + \sum\limits_{A_{ij} < 0} A_{ij} \overline{u}_j$
    \For{$k \in \mathcal{N}$}
      \If{$A_{ik} > 0$} \Comment{Compute implied upper bound.} 
        \State $\overline{u}_k = \min(\overline{u}_k, \overline{l}_k + \frac{b_i - \text{act}_-}{A_{ik}})$
        \If{$k\in\mathcal{I}$}
            $\overline{l}_k = \floor{\overline{l}_k}$
        \EndIf
      \ElsIf{$A_{ik} < 0$} \Comment{Compute implied lower bound.}
        \State $\overline{l}_k = \max(\overline{l}_k, \overline{u}_k + \frac{b_i - \text{act}_-}{A_{ik}})$
        \If{$k\in\mathcal{I}$}
            $\overline{l}_k = \ceil{\overline{l}_k}$
        \EndIf
      \EndIf
      \If{$\overline{l} > \overline{u}$} \Comment{Variable domain became empty.}
        \State \Return (INFEASIBLE, $\overline{l}$, $\overline{u}$)
     \EndIf
      \If{$\overline{u}_k$ or $\overline{l}_k$ changed}
        $Q = Q\cup \{k\}$
      \EndIf
    \EndFor
      \If{$\min\limits_{\overline{l}\leq x \leq \overline{u}}\ A_{il} x > b_i$} \Comment{Row became infeasible.}
        \State \Return (INFEASIBLE, $\overline{l}$, $\overline{u}$)
      \EndIf
  \EndFor
\EndWhile
\State \Return (FEASIBLE, $\overline{l}$, $\overline{u}$)
\end{algorithmic}
\end{algorithm}

\subsubsection{Numerical results}

We implemented the algorithm as a heuristic in the open-source solver \SCIP 9.1.1~\cite{SCIP9} using the Python interface \PySCIPOpt~\cite{PYSCIPOPT}.
This allows easy integration of our pre-trained \NN.
Additionally, \SCIP allows the user to enter the so-called \textit{probing mode}. Probing mode enables quick implementation of a propagation scheme by providing out-of-the-box backtracking and propagation (\cref{alg:propagation}) operators.
Note that commercial solvers do not expose such features.
\GUROBI was used as the \LP solver within \SCIP to counter performance issues posed by the large-scale \LPs that would still need solving after running the heuristic.

While testing the algorithm, we figured out two tweaks that vastly improved performance:
The first tweak is to analyze the \NN output, then sort the variables according to which integer variables are correctly predicted more often.
This prevents premature poor fixings, thereby avoiding infeasibility at that stage. It is preferable to have these variables fixed passively through propagation.
The second trick is to set a higher number of propagation rounds performed by \SCIP, thus increasing the chance of fixing the variables before an incorrect rounding instruction is encountered.

Our approach delivers solid results on test instances not seen during training, as summarized in \cref{table:mlheurresults}.
As indicated by the average dual gap ($\le 0.19 \%$ overall), our heuristic effectively mirrors the performance of commercial \MIP solvers with high accuracy for X-Small and Small instances and reasonably well for Large instances.
While robustness decreases with the size of the instances, a larger number of propagation rounds (up to a sweet spot) usually enhance the quality of the upper bounds and increase the robustness of the heuristic.

\begin{table}[h!]
\small
\begin{center}
\begin{tabular}{l c c c c}
\toprule
Instance set & \# propagation  & Success rate & average gap (\%) & time (s) \\
    &  rounds & (\%) & (solved) & (solved) \\ 
\midrule
X-Small &                     1 &     92 \tablefootnote{The remaining 8 instances were solved directly by \SCIP using the \LP with barrier solution (without crossover).}
                                            &                 0.00 &  4.16 \\ 
\midrule
Small   &                     1 &      11 &                 0.19 &   158.46 \\
       &                     5 &      71 &                 0.14 &   164.05 \\
       &                    10 &      88 &                 0.14 &   164.79 \\
       &                    15 &      96 &                 0.14 &   163.65 \\
\midrule
Large   &                     1 &       8 &                 0.09 & 13,947.61 \\
       &                     5 &      27 &                 0.15 & 13,700.93 \\
       &                    10 &      30 &                 0.15 & 13,700.19 \\
       &                    15 &      28 &                 0.13 & 14,268.58 \\
\bottomrule
\end{tabular}
\caption{Results for the machine learning heuristic on the X-Small, Small and Large instance sets. The shifted geometric mean of the duality gap (shift 1\%) and runtime (shift 1s) is only computed for the solved instances.}
\label{table:mlheurresults}
\end{center}
\end{table}

Although promising, the machine learning approach has several drawbacks.
While data generation is manageable for X-Small and Small instances, it becomes problematic for Large instances: \GUROBI required 9 hours on average to compute a \MIP solution (see \cref{table:fp-results-multiple}) for Large instances and failed at converging within 24 hours for 14 of the 100 instances. This highlights the computational effort to train an accurate classifier for large models.
Therefore we refrained from repeating our experiments on the Medium instance set, as we expect the results to be somewhere between those for the Small and Large sets.
For the X-Large instances, this training overhead proved to be intractable.
Although our heuristic is endowed with range propagation, irrecoverable fixing errors are occasionally encountered. 
Lastly, the approach relies on solving two \LPs, which for larger instances becomes far too expensive.

\subsection{Fix-and-Propagate heuristic}
\label{sec:fix-and-propagate}

So far, our experience is that solving the \LP relaxation of our models often proves prohibitively costly.
A successful heuristic cannot rely on the \LP solution to guide its search.
Instead, we focus on the class of root \LP-free primal \MIP heuristics.
Another critical realization is that propagation plays a vital role in deriving proper fixings (\cref{sec:learning-to-round}) and that many propagation rounds improve the success rate.
Both implications point toward the class of Fix-and-Propagate (\FP) heuristics.
In this section, we devise a \FP heuristic similar to the approaches described in \cite{Achterberg2013, Salvagnin2024_AFixPropagateRepairHeuristicsForMIP, FAPStructure_Gamrath2019, BertholdSAP2015}.

FP heuristics implement a two-stage approach.
First, they iteratively fix the integer variables (following a given strategy) and trigger a propagation operator.
Second, all integer variables are fixed in the original \MIP, and the resulting smaller \LP is solved to achieve a high-quality, feasible solution.

\subsubsection{Description of the generic heuristic}

Our implementation is an extension of the Fix-Propagate-Repair code \cite{Salvagnin2024_AFixPropagateRepairHeuristicsForMIP}.
We did not use the repair mechanism since finding good solutions was more challenging than finding feasible solutions.
The simplified scheme of our algorithm is outlined in \cref{alg:fix-and-propagate}. 
It closely resembles the one in \cite{Salvagnin2024_AFixPropagateRepairHeuristicsForMIP}, which we refer to for more details on the implementation of propagation and the branching scheme.
After presolving the initial \MIP with a commercial \MIP solvere generate a sorted list $\mathcal{I}^\circ$ of integer (and binary) variables according to a score, following a \textit{selection strategy}.
We then start our branching procedure at the root node $(l, u)$.
A node is uniquely defined by a pair of lower and upper bounds.
Keeping track of a stack $S$ of nodes, we first apply propagation to each branching node.
If the node is infeasible, we drop it and continue with the next node from the stack.
Otherwise, we select the following unfixed integer variable in the node according to our order $\mathcal{I}^\circ$.
We then generate two child nodes: one with the selected integer fixed to the bound suggested by our \textit{fixing strategy}, the other fixing it to the opposite bound.
Both child nodes get pushed onto the stack of open nodes, and we continue our procedure.
This depth-first-search branching procedure automatically incorporates \textit{backtracking}.
Should both generated child nodes for a given node prove infeasible, we will automatically revert to the node (by using the stack) and try a different fixing.

In our experiments, we allow for an arbitrary amount of backtracking.
However, we limit the number of fixes our heuristic can perform and the number of infeasible fixes it can encounter during its search.
The branching loop terminates when we find an integer feasible partial solution, hit a node limit, encounter too many infeasible nodes, or no nodes are left in the stack.

\begin{algorithm}[h!]
\caption{Fix-and-Propagate heuristic}
\label{alg:fix-and-propagate}
\begin{algorithmic}
\Require $\MIPProb$
\Ensure \MIP feasible solution or NULL if no solution was found
\State Presolve \MIP
\State Sort integers: $\mathcal{I}^\circ \gets $ selection\_strategy($c$, $A$, $b$, $l$, $u$, $\mathcal{I}$)
\State $S \gets (l, u)$ \Comment{Initialized stack with root-node}
\While{$S \neq \emptyset$ and limits not reached}
  \State $(\hat{l}, \hat{u}) \gets \textbf{Pop}(S)$
  \State (status, $\overline{l}$, $\overline{u})$ $\gets$ propagate($A$, $b$, $l$, $u$, $\mathcal{I}$, $\hat{l}$, $\hat{u}$) \Comment{\Cref{alg:propagation}.}
  \If{status is \allcaps{INFEASIBLE}}
    \State \textbf{continue}
  \EndIf
  \State $\hat{\mathcal{I}}^\circ \gets$ pick indices $j \in \mathcal{I}^\circ$ such that $\overline{l}_j < \overline{u}_j$
  \If{$\hat{\mathcal{I}}^\circ = \emptyset$}
    \State \textbf{break}
  \EndIf
  \State Set $i$ as the first element of $\hat{\mathcal{I}}^\circ$
  \State direction $\gets$ fixing\_strategy($c$, $A$, $b$, $l$, $u$, $\mathcal{I}$, $z_i$)
  \State Generate the child nodes $(\overline{l}^\downarrow, \overline{u}^\downarrow)$, $(\overline{l}^\uparrow, \overline{u}^\uparrow)$
  \begin{equation*}
    (\overline{l}_i^\uparrow, \overline{u}_i^\uparrow) \gets
    \begin{cases}
      (\overline{l}_j, \overline{u}_j) \\
      (\overline{u}_j, \overline{u}_j) 
    \end{cases}
    \quad (\overline{l}_i^\downarrow, \overline{u}_i^\downarrow) \gets
    \begin{cases}
      (\overline{l}_j, \overline{u}_j) & \quad j \neq i \\
      (\overline{l}_j, \overline{l}_j) & \quad j = i
    \end{cases}
  \end{equation*}
  \vspace{1cm}
  \If{direction is \allcaps{LOWER}}
    \State\textbf{Push}($S$, $(\overline{l}^\uparrow, \overline{u}^\uparrow)$); \textbf{Push}($S$, $(\overline{l}^\downarrow, \overline{u}^\downarrow)$)
  \Else
    \State\textbf{Push}($S$, $(\overline{l}^\downarrow, \overline{u}^\downarrow)$); \textbf{Push}($S$, $(\overline{l}^\uparrow, \overline{u}^\uparrow)$)
  \EndIf
\EndWhile
\State Solve resulting \LP: $x_{\text{\MIP}} \gets \LPProb$
\State Postsolve $x_{\text{\MIP}}$
\State \Return $x_{\text{\MIP}}$
\end{algorithmic}
\end{algorithm}
Upon the success of the fix-and-propagate phase, we fix all integer variables to their singleton bounds and solve the resulting smaller \LP with a commercial \LP solver; in theory, it is much easier to solve than the original \MIP or the \LP relaxation.
If this smaller \LP is feasible, its solution is used to create a solution for the original presolved \MIP, which is then postsolved by the \MIP solver.
The presolve step of the commercial solver can further exploit the variables with singleton bounds and relaxed integrality conditions.

\subsubsection{A custom fixing and selection strategy}

All the combinations of fixing and selection strategies described in \cite{Salvagnin2024_AFixPropagateRepairHeuristicsForMIP} are readily available in the original FPR code.
However, our experiments with the \LP-free strategies were not successful.
Fixing variables greedily toward improving objective seemed the most promising variant.
Still, it failed due to the low objective density -- the number of nonzero coefficients in the objective -- of the discrete variables (5 \%), compared to the overall objective density of 30\%.
As a remedy, we implemented a custom strategy for both selection and fixing based on the concept of \textit{inferred objective}.
Moving $x$ towards its upper bound increases the value of $y$, which goes against the improving objective direction of $y$ if $c_y > 0$.
Consequently, the inferred objective of $x$ is set to $c_x + c_y$; both variables contribute.
If, on the other hand, $c_y < 0$, $y$ is moved along its improving objective direction, thus the inferred objective of $x$ is $c_x$.

This strategy deduces an inferred objective for every variable by examining all the rows it appears in (\cref{alg:inferred_objective}).
Equalities are treated as two inequalities.
As the inferred objective propagates through the problem via the rows, multiple updating rounds are carried out until no objective moves away from zero.
An upper bound on the number of updating rounds prevents cycling.
A variable with a zero inferred objective is randomly fixed to one of its bounds.

\begin{algorithm}[h!]
\caption{Inferred objective computation}
\label{alg:inferred_objective}
\begin{algorithmic}
\Require $\MIPProb$
\Ensure Inferred objectives $\hat{c}$
\State $l = 0$
\State Initialize $\hat{c}^l_i = \text{sign}(c_i)$ \Comment{$\text{sign}(0) := 0$.}
\State run $\gets$ true
\While{run}
  \State run $\gets$ false
  \State $\hat{c}^{l+1} \gets \hat{c}^{l}$
  \For{$j=1,\dots,n$}
    \ForAll{$i\in\{0,\dots,m\}$ with $a_{ij} x_j \leq b_i - \sum_{k \neq j} a_{ik} x_k$, $a_{ij} \neq 0$}
        \If {$a_{ij}  \hat{c}^{l}_j > 0$ or $a_{ij}  \hat{c}^{l}_j < 0 $}
            \textbf{continue}
        \EndIf
        \For {$a_{ik} \neq 0$, $k \neq j$}
            \If {$a_{ij}  \hat{c}^{l}_j \leq 0$}
              \textbf{continue}
            \EndIf
            \If {$c^{l+1}_j = 0$}
                run $\gets$ true
            \EndIf
            \State $\hat{c}^{l + 1}_j = \hat{c}^{l + 1}_j + \hat{c}^{l}_k$
        \EndFor
    \EndFor  
  \EndFor
  \State $l=l+1$
\EndWhile
\State \Return $\hat{c}^{l-1}$
\end{algorithmic}
\end{algorithm}

\subsubsection{Numerical results}

We compared our \FP heuristic against the commercial solver \GUROBI on the Small, Medium, Large, and X-Large instances.
The experiments were conducted on a cluster of 32 machines, each equipped with 4 Intel(R) Xeon(R) Gold 6342 CPUs running at 2.8 GHz with 500 GB of RAM.
For all runs, we set a time limit of 2 days, a thread limit of 32, and a memory limit of 200 GB.
\CPLEX 12.10 \cite{CPLEX12} was used for pre- and postsolving the \MIP.
We solved the final (smaller) \LP with one the commercial solvers \COPT 7.1.3 \cite{COPT71} (\IPM without crossover), \XPRESS 9.3 \cite{FicoXpress9}, \GUROBI 11.0 \cite{Gurobi11} (\IPM and 1\% gap), and \CPLEX 12.10, whichever was the fastest for our instances.
Thanks to the gap limit, \GUROBI spends most of its time in the \IPM and the crossover (instead of generating cutting planes and running primal heuristics).

The results of \FP and \GUROBI are gathered in \cref{table:fp-results-multiple}.
The gaps have been computed by taking the best-known dual bound: this is either the final \MIP dual bound found by \GUROBI, or, if \GUROBI could not solve the \MIP, the optimum of the \LP relaxation using the \IPM without crossover.
Across the four instance sets, \FP is significantly faster than \GUROBI, a trend that grows with the instance size.
Furthermore, the success rate of \FP is 100\%, while \GUROBI fails to solve some of the Large and all of the X-Large instances within the time limit (it usually gets stuck in the crossover).
\FP finds near-optimal primal solutions for the Small and Large instances, however the duality gap for the Medium and X-Large instances lies between 20 and 35\%.

\begin{table}[h!]
\small
\begin{center}
\begin{tabular}{llccrr}
\toprule
Instance set & Solver & Solved instances & Success rate (\%) & Avg. time [s] & Avg. gap (\%) \\
\midrule
Small    & \FP     & 1,000 & 100 &    4.01 &  5.68 \\
         & \GUROBI & 1,000 & 100 &   16.90 &  0.41 \\
\midrule
Medium   & \FP     &   96 & 100 &   13.04 & 22.61 \\
         & \GUROBI &   96 & 100 &   77.09 &  0.05 \\
\midrule
Large    & \FP     &  100 & 100 &  333.92 &  1.03 \\
         & \GUROBI &   86 &  86 & 4,548.24 &  0.05 \\
\midrule
X-Large  & \FP     &   20 & 100 & 3,487.27 & 34.96 \\
         & \GUROBI &    0 &   0 &       - &     - \\
\bottomrule
\end{tabular}
\caption{Comparison of \FP and \GUROBI on the Small, Medium, Large and X-Large instance sets. Gap measured to best known dual bound.}
\label{table:fp-results-multiple}
\end{center}
\end{table}

Given its simplicity, the quality of the upper bounds found by the \FP heuristic is remarkable.
It often produces near-optimal solutions in a short amount of time.
It scales much better than the more exhaustive search implemented in commercial solvers since it does not rely on solving the \LP relaxation.

We also had \GUROBI mimic \FP by disabling cuts and crossover and by setting the node limit to 1.
We hoped that \GUROBI would skip the crossover and cutting plane generation and focus on running primal heuristics.
While \GUROBI no longer got stuck in the crossover, it still could not produce any feasible solution within the time limit.

\subsection{Summary}
\begin{figure}[htbp!]
\centering
\begin{minipage}{0.48\textwidth}
    \centering
    \begin{tikzpicture}[scale=0.9]
        \begin{axis}[
            ylabel={Runtime (s)},
            boxplot/draw direction=y,
            xtick={3, 8.0, 13, 18},
            xticklabels={Small, Medium, Large, X-Large},
            extra x ticks={0.5, 5.5, 10.5, 15.5, 20.5}, 
            extra x tick labels={}, 
            extra x tick style={
                tick label style={draw=none}, 
                grid=major, 
                major tick length=0pt, 
            },
            cycle list={{fill=cp1!40,draw=cp1}, {fill=cp2!40,draw=cp2}, {fill=cp3!40,draw=cp3}, {fill=cp5!40,draw=cp5}, {fill=cp4!40,draw=cp4}},
            legend style={draw=none}, 
            ymode=log,
            xtick align=center, 
            axis line style={-}, 
            tick style={draw=none}, 
            grid style={dashed, gray!30},
        ]

\addplot+ [boxplot prepared={median=35.25, upper quartile=40.35, lower quartile=32.91, upper whisker=102.98, lower whisker=23.72}] coordinates {};
\addplot+ [boxplot prepared={median=4.16, upper quartile=4.20, lower quartile=4.10, upper whisker=3.95, lower whisker=4.3}] coordinates {};
\addplot+ [boxplot prepared={median=3.58, upper quartile=4.56, lower quartile=3.12, upper whisker=6.72, lower whisker=2.64}] coordinates {};
\addplot+ [boxplot prepared={median=15.29, upper quartile=15.45, lower quartile=15.11, upper whisker=15.68, lower whisker=15.00}] coordinates {};
\addplot+ [boxplot prepared={median=26.11, upper quartile=31.52, lower quartile=23.03, upper whisker=1431.53, lower whisker=16.49}] coordinates {};

\addplot+ [boxplot prepared={median=131.46, upper quartile=155.74, lower quartile=111.99, upper whisker=389.73, lower whisker=80.29}] coordinates {};
\addplot+ [boxplot prepared={median=162.52, upper quartile=170.39, lower quartile=154.41, upper whisker=194.36, lower whisker=141.23}, fill=white, draw=white] coordinates {}; 
\addplot+ [boxplot prepared={median=13.22, upper quartile=14.36, lower quartile=11.70, upper whisker=17.30, lower whisker=8.39}] coordinates {}; 
\addplot+ [boxplot prepared={median=76.28, upper quartile=85.26, lower quartile=68.61, upper whisker=108.18, lower whisker=45.06}] coordinates {}; 
\addplot+ [boxplot prepared={median=711.51, upper quartile=1099.53, lower quartile=512.98, upper whisker=2880.95, lower whisker=248.60}] coordinates {};

\addplot+ [boxplot prepared={median=24002.62, upper quartile=30818.48, lower quartile=18893.40, upper whisker=126893.06, lower whisker=13611.37}] coordinates {}; 
\addplot+ [boxplot prepared={median=13484, upper quartile=14854, lower quartile=12090, upper whisker=16217, lower whisker=12090}] coordinates {}; 
\addplot+ [boxplot prepared={median=292.56, upper quartile=410.77, lower quartile=234.19, upper whisker=659.21, lower whisker=181.12}] coordinates {}; 
\addplot+ [boxplot prepared={median=2025.13, upper quartile=3839.62, lower quartile=1378.03, upper whisker=7517.62, lower whisker=975.67}] coordinates {}; 
\addplot+ [boxplot prepared={median=23573.86, upper quartile=172800, lower quartile=15422.63, upper whisker=172800, lower whisker=9704.30}] coordinates {};

\addplot+ [boxplot prepared={median=15, upper quartile=22, lower quartile=12, upper whisker=28, lower whisker=10}, fill=white, draw=white] coordinates {}; 
\addplot+ [boxplot prepared={median=172800, upper quartile=172800, lower quartile=172800, upper whisker=172800, lower whisker=172800}, fill=white, draw=white] coordinates {};
\addplot+ [boxplot prepared={median=3456.78, upper quartile=3671.04, lower quartile=3157.15, upper whisker=4441.89, lower whisker=2659.68}] coordinates {}; 
\addplot+ [boxplot prepared={median=172800, upper quartile=202800, lower quartile=172800, upper whisker=172800, lower whisker=172800}, fill=white, draw=white] coordinates {};
\addplot+ [boxplot prepared={median=172800, upper quartile=172800, lower quartile=172800, upper whisker=172800, lower whisker=172800}, fill=white, draw=white] coordinates {};

        \end{axis}
    \end{tikzpicture}
\end{minipage}\hfill
\begin{minipage}{0.48\textwidth}
    \centering
    \begin{tikzpicture}[scale=0.9]
        \begin{axis}[
            ylabel={Gap in \%},
            boxplot/draw direction=y,
            xtick={3, 8.0, 13, 18},
            xticklabels={Small, Medium, Large, X-Large},
            cycle list={{fill=cp1!40,draw=cp1}, {fill=cp2!40,draw=cp2}, {fill=cp3!40,draw=cp3}, {fill=cp5!40,draw=cp5}, {fill=cp4!40,draw=cp4}},
            legend style={draw=none}, 
            extra x ticks={0.5, 5.5, 10.5, 15.5, 20.5}, 
            extra x tick labels={}, 
            extra x tick style={
                tick label style={draw=none}, 
                grid=major, 
                major tick length=0pt, 
            },
            xtick align=center, 
            axis line style={-}, 
            ymode=log,
            ymin=1e-5,
            ymax=100,
            tick style={draw=none}, 
            grid style={dashed, gray!30},
        ]
\addplot+ [boxplot prepared={median=3.41e-05, upper quartile=5.28e-05, lower quartile=1.68e-5, upper whisker=9.97e-05, lower whisker=1.18e-8}] coordinates {};
\addplot+ [boxplot prepared={median=0.002, upper quartile=0.005, lower quartile=0.001, upper whisker=0.01, lower whisker=0.00001}] coordinates {};
\addplot+ [boxplot prepared={median=5.68, upper quartile=6.26, lower quartile=5.10, upper whisker=8.00, lower whisker=3.36}] coordinates {};
\addplot+ [boxplot prepared={median=0.37, upper quartile=0.53, lower quartile=0.29, upper whisker=0.88, lower whisker=0.001}] coordinates {};
\addplot+ [boxplot prepared={median=3.20e-5, upper quartile=4.59e-5, lower quartile=2.13e-5, upper whisker=9.99e-5, lower whisker=1.03e-7}] coordinates {};

\addplot+ [boxplot prepared={median=3.41e-5, upper quartile=4.84e-5, lower quartile=2.70e-6, upper whisker=9.87e-5, lower whisker=2.41e-7}] coordinates {}; 
\addplot+ [boxplot prepared={median=0.05, upper quartile=0.06, lower quartile=0.03, upper whisker=0.09, lower whisker=0.0001}, draw=white, fill=white] coordinates {}; 
\addplot+ [boxplot prepared={median=22.17, upper quartile=26.1, lower quartile=18.76, upper whisker=37.01, lower whisker=12.24}] coordinates {}; 
\addplot+ [boxplot prepared={median=0.05, upper quartile=0.06, lower quartile=0.03, upper whisker=0.09, lower whisker=0.0001}] coordinates {}; 
\addplot+ [boxplot prepared={median=2.19e-5, upper quartile=4.83e-5, lower quartile=6.19e-6, upper whisker=0.0001, lower whisker=2.87e-6}] coordinates {};

\addplot+ [boxplot prepared={median=1.39e-5, upper quartile=3.43e-5, lower quartile=3.11e-6, upper whisker=0.90, lower whisker=5.97e-7}] coordinates {}; 
\addplot+ [boxplot prepared={median=0.17, upper quartile=0.19, lower quartile=0.05, upper whisker=0.007, lower whisker=0.22}] coordinates {};
\addplot+ [boxplot prepared={median=0.7, upper quartile=1.3, lower quartile=0.48, upper whisker=2.59, lower whisker=0.16}] coordinates {}; 
\addplot+ [boxplot prepared={median=0.01, upper quartile=0.7, lower quartile=0.001, upper whisker=1.6, lower whisker=0.0001}] coordinates {}; 
\addplot+ [boxplot prepared={median=1.50e-5, upper quartile=5.28e-5, lower quartile=4.50e-6, upper whisker=0.96, lower whisker=3.13e-7}] coordinates {};

\addplot+ [boxplot prepared={median=15, upper quartile=22, lower quartile=12, upper whisker=28, lower whisker=10}, fill=white, draw=white] coordinates {}; 
\addplot+ [boxplot prepared={median=0.05, upper quartile=0.06, lower quartile=0.03, upper whisker=0.09, lower whisker=0.0001}, draw=white, fill=white] coordinates {}; 
\addplot+ [boxplot prepared={median=24.99, upper quartile=32.52, lower quartile=21.05, upper whisker=47.23, lower whisker=16.69}] coordinates {}; 
\addplot+ [boxplot prepared={median=0.05, upper quartile=0.06, lower quartile=0.03, upper whisker=0.09, lower whisker=0.0001}, draw=white, fill=white] coordinates {};
\addplot+ [boxplot prepared={median=0.05, upper quartile=0.06, lower quartile=0.03, upper whisker=0.09, lower whisker=0.0001}, draw=white, fill=white] coordinates {}; 

        \end{axis}
    \end{tikzpicture}
\end{minipage}
\begin{tikzpicture}
    \begin{axis}[
        hide axis,
        xmin=0, xmax=1, ymin=0, ymax=1,
        legend columns=5,
        legend style={at={(0.5,-0.2)}, anchor=north, draw=black, 
                      rounded corners=2pt, legend image post style={line width=1.5pt}}, 
    ]
        \addlegendimage{cp1}
        \addlegendentry{\RENS}

        \addlegendimage{cp2}
        \addlegendentry{Machine learning}

        \addlegendimage{cp3}
        \addlegendentry{\FP (\COPT)}

        \addlegendimage{cp5}
        \addlegendentry{\GUROBI 1\%}

        \addlegendimage{cp4}
        \addlegendentry{\GUROBI}
    \end{axis}
\end{tikzpicture}
\caption{Solving time and gap for all strategies and solved instances over the Small to X-Large sets.}
\label{fig:final_results}
\end{figure}
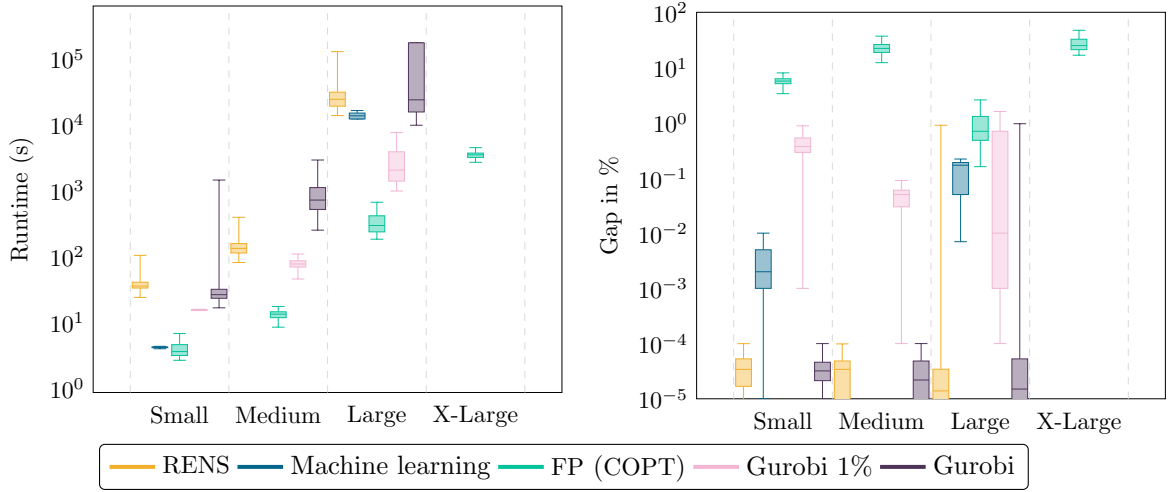

 gather the results of all three methods on the Small, Medium, Large, and X-Large instance sets in \cref{fig:final_results}.
Both box plots have a log-scaled y-axis.
The left plot shows runtimes and the right plot displays the duality gap at the solution.
Some results were not computed (e.g., machine learning on the Medium set and \GUROBI on the X-Large set); we left the corresponding columns empty.

Devising the \FP strategy allowed us to relabel the Large instance set from "possible to solve" to "practically applicable".
Furthermore, we were able to produce feasible solutions for all X-Large instances, to which all commercial solvers failed.

\section{Conclusion and perspectives}

This paper presented three alternative solution strategies for approximately solving \ESOM \MIP instances that were, up to now, intractable for real-world applications.
Thanks to a preliminary analysis of our models, we moved away from standard branch-and-bound \MIP solution techniques and opted for primal heuristics.
We incrementally refined our developments -- with varying success: i) the \RENS heuristic exploits the solution to the \LP relaxation and defines a sub-\MIP with a significantly reduced search space; ii) our machine learning strategy combines a neural fixing strategy and range propagation; iii) a novel \LP-free \FP heuristic, based on inferred costs, also exploits the strength of range propagation.
We demonstrated how to effectively tackle intractable instances: our \FP heuristic solves our largest instances (82M nonzeros), often to acceptably small duality gaps and always within a reasonable time.

We wish to emphasize two points.
First, we believe that any model instance with similar properties to our instances should be solvable (with a few modifications) with the \LP-free \FP heuristic or a \RENS neighborhood approach.
Second, classical \MIP solvers might not be the correct tool for solving the ever-growing class of \ESOMs when high accuracy is not paramount.
Instead, we demonstrated that tailored heuristics can quickly produce high-quality solutions.
This should enable decision-makers and modelers to integrate larger and more scenarios into their models to mitigate data uncertainty, thus enhancing their models' applicability.

We hope practitioners can benefit from our detailed development process and discussion.
We intend to pursue our research on primal heuristics tailored to \ESOM \MIPs and build a portfolio of primal techniques for the quick solution of \ESOMs.
An avenue for reflection is to determine which problem features are required for a neural fixing strategy to perform well without the prior knowledge of the \LP solution.
We also plan to experiment with more sophisticated \NN architectures and learning methods that may scale more effectively.
The similar underlying structure of the UNSEEN instances could be explored with Graph Convolutional \NNs (\GCNNs)~\cite{Zhang2019}.
\GCNNs have proven successful in solving \MIPs by learning branching~\cite{Gasse2019} and cut selection~\cite{Turner2022}.
Unlike multi-layer perceptrons, \GCNNs allow input graphs of variable size: a \GCNN can be trained on smaller instances and used for prediction on larger instances for which generating training data is harder and more expensive.
This most likely has limitations since generalization on larger instances is not simple.

\section*{Acknowledgements}
{\footnotesize The described research activities are funded by the Federal Ministry for Economic Affairs and Energy within the project UNSEEN (ID: 03EI1004A, 03EI1004D).
Part of the work for this article has been conducted in the Research Campus MODAL funded by the German Federal Ministry of Education and Research (BMBF) (fund numbers 05M14ZAM, 05M20ZBM). The work has been co-funded by the European Union (European Regional Development Fund ERDF, fund number: STIIV-001).

The authors thank Domenico Salvagnin for making his fix-propagate-repair code available.}

\bibliographystyle{plain}
\bibliography{references}

\end{document}